\documentclass[]{article}

\usepackage[a4paper, portrait, margin=1.1811in]{geometry}
\usepackage[english]{babel}
\usepackage[utf8]{inputenc}
\usepackage[T1]{fontenc}
\usepackage{helvet}
\usepackage{etoolbox}
\usepackage{graphicx}
\usepackage{titlesec}
\usepackage[font=footnotesize]{caption}
\usepackage{booktabs}
\usepackage{xcolor} 
\usepackage[colorlinks, citecolor=teal, linkcolor=teal, urlcolor=teal]{hyperref}
\usepackage{caption}
\graphicspath{ {./images/} }
\usepackage{scrextend}
\usepackage{fancyhdr}
\usepackage{authblk}
\usepackage{multicol}
\usepackage{comment}
\usepackage{float}

\usepackage{multirow}
\usepackage{subfig}
\usepackage{makecell} 

\usepackage{lmodern}

\usepackage{amsmath}
\usepackage{amssymb}

\usepackage{amsthm}

\usepackage{xcolor}
\usepackage[linesnumbered,ruled,vlined]{algorithm2e}

\usepackage{makeidx}
\usepackage[acronym, nonumberlist, style=long, toc]{glossaries}

\usepackage{xkeyval, xfor, amsgen, array}

\usepackage{natbib}   
\setcitestyle{authoryear, open={(},close={)}}

\usepackage{orcidlink}
\usepackage[misc]{ifsym}

\usepackage{fontawesome}

\usepackage{soul}
\usepackage[framemethod=TikZ]{mdframed}

\theoremstyle{definition}

\SetCommentSty{mycommfont}

\SetKwInput{KwInput}{Input}                
\SetKwInput{KwOutput}{Output}              

\newcounter{defi}[section]
\renewcommand{\thedefi}{\arabic{section}.\arabic{defi}}
\newenvironment{defi}[2][]{%
\refstepcounter{defi}%
\ifstrempty{#1}%
{\mdfsetup{%
frametitle={%
\tikz[baseline=(current bounding box.east),outer sep=0pt]
\node[anchor=east,rectangle,fill=gray!20]
{\strut Definition~\thedefi};}}
}%
{\mdfsetup{%
frametitle={%
\tikz[baseline=(current bounding box.east),outer sep=0pt]
\node[anchor=east,rectangle,fill=gray!20]
{\strut Definition~\thedefi:~\textbf{#1}};}}%
}%
\mdfsetup{innertopmargin=10pt,linecolor=gray!20,%
linewidth=2pt,topline=true,%
frametitleaboveskip=\dimexpr-\ht\strutbox\relax
}
\begin{mdframed}[]\relax%
\label{#2}}{\end{mdframed}}

\newcounter{theo}[section]
\renewcommand{\thetheo}{\arabic{section}.\arabic{theo}}
\newenvironment{theo}[2][]{%
\refstepcounter{theo}%
\ifstrempty{#1}%
{\mdfsetup{%
frametitle={%
\tikz[baseline=(current bounding box.east),outer sep=0pt]
\node[anchor=east,rectangle,fill=gray!20]
{\strut Theorem~\thetheo};}}
}%
{\mdfsetup{%
frametitle={%
\tikz[baseline=(current bounding box.east),outer sep=0pt]
\node[anchor=east,rectangle,fill=gray!20]
{\strut Theorem~\thetheo:~#1};}}%
}%
\mdfsetup{innertopmargin=10pt,linecolor=gray!20,%
linewidth=2pt,topline=true,%
frametitleaboveskip=\dimexpr-\ht\strutbox\relax
}
\begin{mdframed}[]\relax%
\label{#2}}{\end{mdframed}}

\newcounter{prf}[section]

\newenvironment{prf}[2][]{%
\refstepcounter{prf}%
\ifstrempty{#1}%
{\mdfsetup{%
frametitle={%
\tikz[baseline=(current bounding box.east),outer sep=0pt]
\node[anchor=east,rectangle,fill=gray!20]
{\strut Proof};}}
}%
{\mdfsetup{%
frametitle={%
\tikz[baseline=(current bounding box.east),outer sep=0pt]
\node[anchor=east,rectangle,fill=gray!20]
{\strut Proof:~#1};}}%
}%
\mdfsetup{innertopmargin=10pt,linecolor=gray!20,%
linewidth=2pt,topline=true,%
frametitleaboveskip=\dimexpr-\ht\strutbox\relax
}
\begin{mdframed}[]\relax%
\label{#2}}{\qed\end{mdframed}}

\newcounter{rema}[section]
\renewcommand{\therema}{\arabic{section}.\arabic{rema}}
\newenvironment{rema}[2][]{%
\refstepcounter{rema}%
\ifstrempty{#1}%
{\mdfsetup{%
frametitle={%
\tikz[baseline=(current bounding box.east),outer sep=0pt]
\node[anchor=east,rectangle,fill=gray!20]
{\strut Remark~\therema};}}
}%
{\mdfsetup{%
frametitle={%
\tikz[baseline=(current bounding box.east),outer sep=0pt]
\node[anchor=east,rectangle,fill=gray!20]
{\strut Remark~\therema:~#1};}}%
}%
\mdfsetup{innertopmargin=10pt,linecolor=gray!20,%
linewidth=2pt,topline=true,%
frametitleaboveskip=\dimexpr-\ht\strutbox\relax
}
\begin{mdframed}[]\relax%
\label{#2}}{\end{mdframed}}

\newcounter{prop}[section]
\renewcommand{\theprop}{\arabic{section}.\arabic{prop}}

\makeglossaries

\newacronym{BDCA}{BDCA}{Boosted Difference of Convex Functions Algorithm}

\newacronym[longplural={Confidence Intervals}, \glsshortpluralkey={CIs}]{CI}{CI}{Confidence Interval}

\newacronym[longplural={Conditional Value-at-Risks}, \glsshortpluralkey={CVaRs}]{CVaR}{CVaR}{Conditional Value-at-Risk}

\newacronym{CDF}{CDF}{Cumulative Distribution Function}

\newacronym{DC}{DC}{Difference of Convex Functions}
\newacronym{DCA}{DCA}{Difference of Convex Functions Algorithm}

\newacronym{dkabs}{$\mathbf{d}_k$-abs}{Absolute $\Delta$ in Solution Vector}
\newacronym{dkrel}{$\mathbf{d}_k$-rel}{Relative $\Delta$ in Solution Vector}

\newacronym{DJ}{DJ}{Dow Jones Industrial Average}

\newacronym{FTSE}{FTSE}{Financial Times Stock Exchange 100}
\newacronym{FF49}{FF49}{Fama \& French 49 Industrial Portfolios}

\newacronym{fctabs}{$\phi$-abs}{Absolute $\Delta$ in Objective Function}
\newacronym{fctcrel}{$\phi$-rel}{Relative $\Delta$ in Objective Function}

\newacronym{IQR}{IQR}{Interquartile Range}

\newacronym{KKT}{KKT}{Karush–Kuhn–Tucker}
\newacronym{KL}{KL}{Kurdyka-\L ojasiewicz}

\newacronym{LICQ}{LICQ}{Linear Independence Constraint Qualification}

\newacronym{MFCQ}{MFCQ}{Mangasarian-Fromovitz Constraint Qualification}

\newacronym{MV}{MV}{Mean-Variance}

\newacronym{NASDAQ}{NASDAQ}{National Association of Securities Dealers Automated Quotation 100}

\newacronym{pDCA}{pDCA}{proximal Difference of Convex Functions Algorithms}

\newacronym{PDF}{PDF}{Probability Density Function}

\newacronym{SCA}{SCA}{Sequential Convex Approximations}

\newacronym{SAA}{SAA}{Sample Average Approximation}

\newacronym{VaR}{VaR}{Value-at-Risk}

\begin{document}

\pagestyle{fancy}
\lhead{Thormann et al.}
\rhead{BDCA for VaR Constrained Portfolio Optimization}

\newpage
\setcounter{page}{1}
\renewcommand{\thepage}{\arabic{page}}
\setlength{\headheight}{35.2pt}
	
\captionsetup[figure]{labelfont={bf},labelformat={default},labelsep=period,name={Figure }}	

\captionsetup[table]{labelfont={bf},labelformat={default},labelsep=period,name={Table }}

\setlength{\parskip}{0.5em}



\fancypagestyle{plain}{
	\fancyhf{}
	\renewcommand{\headrulewidth}{0pt}
	\renewcommand{\familydefault}{\sfdefault}
	
	\lhead{\color{cyan}\normalsize \textbf{Working Paper} (v2)\\ \color{black}
	\text{Operational Research Group, University of Southampton}\\ }
	
}

\makeatletter
\patchcmd{\@maketitle}{\LARGE \@title}{\fontsize{19.15}{24}\selectfont\@title}{}{}
\makeatother

\setlength{\affilsep}{1em}  
\newsavebox\affbox

\author{
    {\large \textbf{Marah-Lisanne Thormann}} {\normalsize \href{mailto:m.-l.thormann@soton.ac.uk}{\Letter} \href{https://www.southampton.ac.uk/people/5ztphy/ms-marah-thormann}{\faHome} \protect\\ \vspace{0.1cm}
     School of Mathematical Sciences, \protect\\ University of Southampton, \protect\\ 
     SO17 1BJ Southampton, UK}  \vspace{0.1cm} \\
     
    {\large \textbf{Phan Tu Vuong}} {\normalsize {\small\href{mailto:t.v.phan@soton.ac.uk}{\Letter}}  \href{https://www.southampton.ac.uk/people/5y2ds9/doctor-vuong-phan}{\faHome} {\large \orcidlink{0000-0002-1474-994X}} \protect\\ \vspace{0.1cm}
     School of Mathematical Sciences, \protect\\
     University of Southampton, \protect\\ SO17 1BJ Southampton, UK}  \vspace{0.25cm}
    
    {\large \textbf{Alain B. Zemkoho}} {\normalsize {\small \href{mailto:a.b.zemkoho@soton.ac.uk}{\Letter}} \href{https://www.southampton.ac.uk/~abz1e14/index.html}{\faHome} {\large \orcidlink{0000-0003-1265-4178}} \protect\\
     \vspace{0.1cm} School of Mathematical Sciences, \protect\\
     University of Southampton, \protect\\ SO17 1BJ Southampton, UK}
}

\titlespacing\section{0pt}{12pt plus 4pt minus 2pt}{0pt plus 2pt minus 2pt}

\titlespacing\subsection{12pt}{12pt plus 4pt minus 2pt}{0pt plus 2pt minus 2pt}

\titlespacing\subsubsection{12pt}{12pt plus 4pt minus 2pt}{0pt plus 2pt minus 2pt}

\titleformat{\section}{\normalfont\fontsize{10}{15}\bfseries}{\thesection.}{1em}{}
\titleformat{\subsection}{\normalfont\fontsize{10}{15}\bfseries}{\thesubsection.}{1em}{}
\titleformat{\subsubsection}{\normalfont\fontsize{10}{15}\bfseries}{\thesubsubsection.}{1em}{}


\title{\vspace{-0.5cm} \textbf{The Boosted Difference of Convex Functions Algorithm for Value-at-Risk Constrained Portfolio Optimization} \vspace{0.4cm}}
\date{}  

\begingroup
\let\center\flushleft
\let\endcenter\endflushleft

\maketitle

\vspace{-0.75cm}

	

\noindent\rule{1.25cm}{0.4pt}  \rlap{\color{gray}\vrule}%
  \fboxsep1.5mm\colorbox[rgb]{1,1,1}{\raisebox{-0.4ex}{%
    \large\selectfont\sffamily\bfseries\abstractname}} \noindent\rule{11.6cm}{0.4pt}

    \noindent
    A highly relevant problem of modern finance is the design of \gls{VaR} optimal portfolios. 
    Due to contemporary financial regulations, banks and other financial institutions are tied to use the risk measure to control their credit, market, and operational risks. 
    Despite its practical relevance, the non-convexity induced by \gls{VaR} constraints in portfolio optimization problems remains a major challenge.
    To address this complexity more effectively, this paper proposes the use of the \gls{BDCA} to approximately solve a Markowitz-style portfolio selection problem with a \gls{VaR} constraint. 
    As one of the key contributions, we derive a novel line search framework that allows the application of the algorithm to \gls{DC} programs where both components are non-smooth.
    Moreover, we prove that the \gls{BDCA} linearly converges to a Karush-Kuhn-Tucker point for the problem at hand using the Kurdyka– \L ojasiewicz property.
    We also outline that this result can be generalized to a broader class of piecewise-linear \gls{DC} programs with linear equality and inequality constraints. 
    In the practical part, extensive numerical experiments under consideration of best practices then demonstrate the robustness of the \gls{BDCA} under challenging constraint settings and adverse initialization. 
    In particular, the algorithm consistently identifies the highest number of feasible solutions even under the most challenging conditions, while other approaches from chance-constrained programming lead to a complete failure in these settings.
    Due to the open availability of all data sets and code, this paper further provides a practical guide for transparent and easily reproducible comparisons of \gls{VaR}-constrained portfolio selection problems in Python. \\ \\
    \noindent 
    \textbf{\text{Keywords}}: Boosted Difference of Convex Functions Algorithm; Chance-Constrained Programming; Kurdyka– \L ojasiewicz Property; Portfolio Selection; Value-At-Risk
    \vspace{0.2cm}
  
    \noindent 
    \textbf{Mathematics Subject Classification (2020)}: 65K05; 65K10; 90C15; 90C26; 91G10; 91G70

    \vspace*{\fill}

    \noindent \rule{15cm}{0.2pt}
    {\footnotesize \noindent 
    \textbf{Funding Acknowledgement}: 
    The first author receives a PhD studentship from the School of Mathematical Sciences at the University of Southampton.}

\endgroup

\glsresetall


\newpage
\section{Introduction}\label{sec:Intro}

``There's a way to do it better ... find it'' [Thomas A. Edison]. 
In line with this maxim, both practitioners and theoreticians have, for over 70 years, sought to improve the foundational \gls{MV} model proposed by \citet{markowitz1952} to solve the portfolio selection problem [cf. \cite{zhang2018}]. 
Over time, two major research streams have emerged within modern portfolio theory: (i) the use of alternative (non-variance-based) risk measures, and (ii) the inclusion of additional criteria and constraints within the optimization framework [cf. \citet[p. 1285]{ANAGNOSTOPOULOS2010}; \citet[p. 105]{RIGHI2018}].
As each refinement increases model complexity, contemporary research often grapples with the trade-off between practical applicability and the accurate representation of complex and unpredictable financial markets [cf. \citet[pp. 138--139]{zhang2018}]. 
In this paper, we address this challenge and demonstrate how the portfolio selection problem, when formulated using an alternative risk measure, can be solved more robustly under demanding risk thresholds than with existing approaches.


To further contextualize our contribution and motivate the proposed methodology, the remainder of this introduction is structured as follows. 
We begin by revisiting the classical \gls{MV} framework and its limitations. 
Next, we introduce \gls{VaR} as an alternative risk measure to the variance constraint and reformulate the optimization problem as a chance-constrained program. 
This is followed by a discussion of the methodological challenges associated with chance constraints.
Afterwards, we review existing solution approaches and highlight the potential of \gls{DC} programming to address these challenges. 
Ultimately, the introduction concludes with a brief overview of the key contributions of this paper and an outline of the subsequent sections.


\subsection{The Mean-Variance Framework and its Limitations}


In the basic framework of the \gls{MV} portfolio, an investor seeks to allocate capital among $n \in \mathbb{N}$ available assets, assuming that no risk-free alternatives are available [cf. \citet[p. 81]{markowitz1952}]. 
This budget must be allocated such that the resulting portfolio is efficient in terms of both expected return and risk [cf. \citet[p. 82]{markowitz1952}]. 
In other words, an efficient portfolio is one for which no alternative allocation yields a higher expected return for a fixed level of risk, or a lower level of risk for the same expected return. 
More technically, this objective can be modeled by defining the return of each individual asset~$i \in \{1, \ldots, n\}$ as a random variable $\xi_i$, where $\mathbb{E} \left[ \xi_i \right] = \mu_i \in \mathbb{R}$ denotes the expected return of the investment and $\mathbb{V} \left[\xi_i \right] = \sigma_{i}^2 \in \mathbb{R}_{+}$ serves as a measure of risk.
By introducing portfolio weights $\mathbf{w} \in \Delta$, with
\begin{equation}
    \Delta := \Big\{ \mathbf{w} = \begin{bmatrix} w_1, \ldots, w_n \end{bmatrix}^\top \in [0, 1]^n  \ \bigl| \  \textstyle \sum_{i=1}^n w_i = 1 \Big\},
\end{equation}
the overall expected return and variance of the portfolio are given by $\mathbb{E}\left[\mathbf{w}^\top \boldsymbol{\xi}\right] = \sum_{i = 1}^{n} w_i \cdot \mu_i$ and $\mathbb{V} \left[\mathbf{w}^\top \boldsymbol{\xi} \right] = \sum_{i = 1}^{n}\sum_{j = 1}^{n} \rho_{ij} \cdot w_i \cdot w_j \cdot \sigma_{i} \cdot \sigma_j,$ where $\boldsymbol{\xi} := \begin{bmatrix} \xi_1, \ldots, \xi_n\end{bmatrix}^\top$ is the vector of random returns, and $\rho_{ij} \in [-1, 1]$ denotes the correlation coefficient between assets $i$ and $j$. 
The term $\mathbf{w}^\top \boldsymbol{\xi} = \sum_{i=1}^n w_i \cdot \xi_i$ represents the weighted sum of individual returns (i.e., the portfolio return), and each weight $w_i$ specifies the proportion of capital allocated to asset $i$. 
Due to the typically imperfect correlation $\rho_{ij}$ between distinct returns $\xi_i$ and $\xi_j$ for $i \neq j$, with $j \in \{1, \ldots, n\}$, combining multiple assets in a portfolio can reduce overall risk through diversification [cf. \citet[chapter~6]{jiang2022}]. 
This effect can be explained by decomposing the total risk into systematic and specific components, where the latter tend to cancel each other out, especially when assets are selected across different industries [cf. \citet[chapter 6]{jiang2022}; \citet[p. 89]{markowitz1952}]. 
Although systematic risk remains, diversification still lowers portfolio variance while preserving expected return [cf. \citet[chapter 6]{jiang2022}]. 
Accordingly, the economic decision problem under uncertainty can be formulated as a constrained stochastic program of the form
\begin{align}\label{eq:markowitz}
    \begin{split}
        \underset{\mathbf{w} \in \Delta}{\text{maximize}} \quad &\mathbb{E}\left[\mathbf{w}^\top \boldsymbol{\xi} \right] \\
        \text{subject to} \quad &\mathbb{V}\left[\mathbf{w}^\top \boldsymbol{\xi} \right] \leq \sigma_{\text{max}}^2,
    \end{split}
\end{align}
where the objective is to determine the portfolio weights $w_i$ such that, for a given risk tolerance $\sigma_{\text{max}} \in \mathbb{R}_+$, the expected return is maximized. 
Alternatively, one may also equivalently formulate the optimization problem by minimizing portfolio variance for a given level of expected return. 
The set of optimal portfolios that achieve the highest expected return for a given level of risk, or the lowest possible risk for a given expected return, is known in the literature as the \textit{efficient frontier} [cf. \citet[pp. 87--89]{markowitz1952}].


Based on the model formulation presented in \eqref{eq:markowitz}, \citet{markowitz1952} introduced a foundational concept for optimal portfolio selection. 
However, several literature reviews on modern portfolio theory --  such as those by \citet{gunjan2023}, \citet{xidonas2020}, \citet{sun2019}, \citet{zhang2018}, \citet{rather2017}, \citet{MANSINI2014}, and \citet{roman2009} --  highlight limitations of the original approach and reference alternative models proposed by subsequent researchers to better capture the complexity of financial markets. 
Common points of criticism include the restricted planning horizon (single- vs. multi-period), the nature of the objective function (single- vs. multi-objective), assumptions about investor risk preferences (risk-neutral vs. risk-averse), distributional assumptions on asset returns (normal vs. non-normal), the omission of practical constraints such as transaction costs or sparsity, and the choice of risk measure [cf. \citet{gunjan2023}; \citet{sun2019}; \citet{zhang2018}]. 
While most of these limitations relate to the first research stream -- the integration of additional real-world features -- the selection of an appropriate risk measure forms a separate branch of research [cf. \citet[p. 1285]{ANAGNOSTOPOULOS2010}]. 
At the core of this second research stream lies the fundamental critique of variance: it symmetrically penalizes all deviations from the mean, thereby treating both very high and low returns as equally undesirable [cf.  \citet{rather2017}; \citet{Gambrah2014}; \citet{hoe2010}]. 
This symmetric treatment of return deviations then stands in contradiction to the overall objective of the investor - the maximization of the expected portfolio return [cf. \citet[p. 98]{Capinski2014}]. 
To address this limitation, the next subsection introduces a \gls{VaR}-based formulation that serves as risk measure in this paper. 


\subsection{Portfolio Optimization under Value-at-Risk Constraint}

\begin{figure}[H]%
    \centering
    \includegraphics[width=1\textwidth]{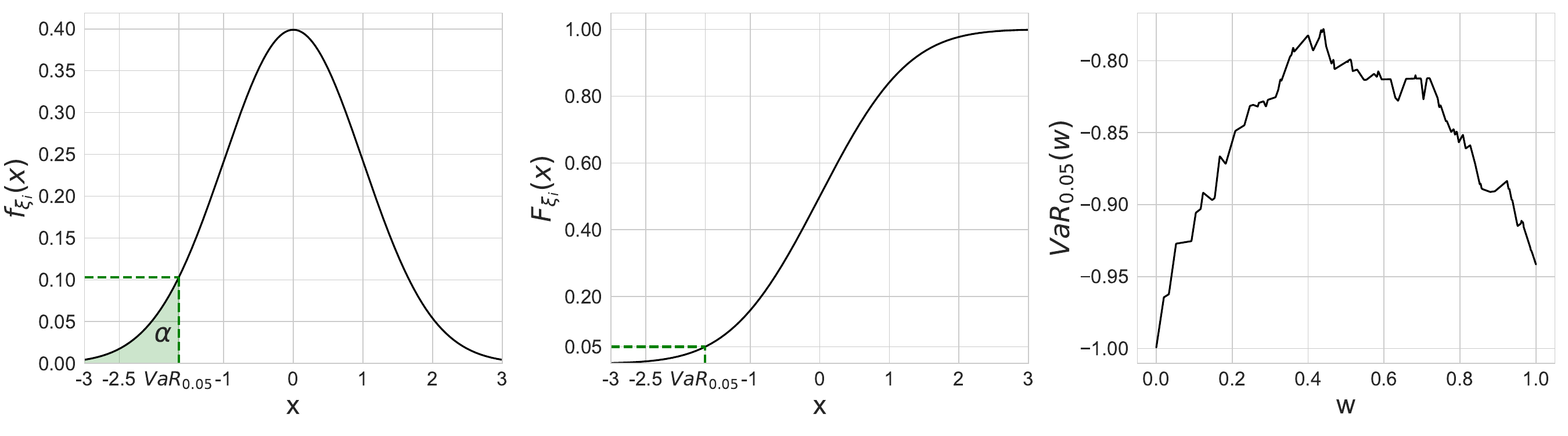}
    \captionsetup{justification=centering}
    \caption{Basic examples designed to explain the concept of $\text{\gls{VaR}}_\alpha$.}%
    \label{fig:var_plots}%
\end{figure}

A prominent alternative to the variance constraint is a corresponding \gls{VaR} formulation, as outlined in the review of portfolio risk measures by \citet{sun2019}. 
Originally introduced by the global financial services company J.P. Morgan Chase~\& Co, \gls{VaR} focuses solely on the downside tail of the return distribution [cf. \citet[p. 1104]{Poncet2022}, \citet[p. 862]{wozabal2012}]. 
Intuitively, the $\text{\gls{VaR}}_{\alpha}$  at confidence level $\alpha \in (0,1)$ represents the lowest return that is not exceeded with probability $1 - \alpha$, assuming unchanged market conditions. 
On a more technical level,  the $\text{\gls{VaR}}_{\alpha}$ of a random return variable $\xi_i$ is defined as the smallest quantile for which the corresponding \gls{CDF} is greater than or equal to~$\alpha$ [cf. \citet[chapter 27]{Poncet2022}]. 
Formally, this can be expressed as
\begin{align}
    \begin{split}
       \text{\gls{VaR}}_\alpha \left(\xi_i \right)  & : = \inf \bigl \{u: F_{\xi_i} \left(u \right) \geq \alpha \bigr \} = F_{\xi_i}^{-1}(\alpha) \quad \text{with} \quad  \alpha =  \int_{- \infty}^{\text{\gls{VaR}}_\alpha} f_{\xi_i} (u) \ \text{d}u, 
    \end{split}
\end{align}
where $f_{\xi_i}: \mathbb{R} \rightarrow \mathbb{R}_+$ denotes the \gls{PDF}, $F_{\xi_i}: \mathbb{R} \rightarrow [0, 1]$ the \gls{CDF}, and $F_{\xi_i}^{-1}: [0, 1] \rightarrow \mathbb{R}$ the inverse \gls{CDF} of $\xi_i$, also referred to as quantile function. 
To graphically illustrate the $\text{\gls{VaR}}_{\alpha}$, Figure~\ref{fig:var_plots} shows the \gls{PDF} (left panel) and \gls{CDF} (middle panel) of a standard normal distribution, where the $\text{\gls{VaR}}_{0.05}$ is indicated by green dashed lines.
Note that in this paper, $\text{\gls{VaR}}_{\alpha}$ is considered to be a (non-concave) acceptability-type functional, but equivalently, it could also be formulated as a (non-convex) deviation-risk-type functional if the random variables represented portfolio losses instead of returns [cf. \citet[subsection 2.2]{Pflug2007}]. 


Beyond its theoretical foundation, the $\text{\gls{VaR}}_{\alpha}$ also offers practical advantages in the context of regulatory compliance. 
Contemporary financial regulations require banks and other financial institutions to maintain minimum levels of capital to ensure solvency in the event of significant losses arising from market, credit, or operational risks [cf. \citet[p. 863]{wozabal2012}]. 
Within this framework, the Basel~II and Solvency~II regulations mandate that the required regulatory capital be determined based on the $\text{\gls{VaR}}_{\alpha}$ of the institution’s key business units [cf. \citet[p. 272]{Sarykalin2014}; \citet[p. 1928]{SANTOS2012}; \citet[p. 362-363]{CUOCO2006}]. 
Consequently, from a risk management perspective, it is both necessary and desirable to actively monitor and control the $\text{\gls{VaR}}_{\alpha}$ of financial portfolios, which further motivates the relevance of the optimization problem considered in this paper [cf. \citet[p. 863]{wozabal2012}].


\vspace{-0.2cm}


In the academic literature, numerous studies - including those by \citet{liu2021}, \citet{MOHAMMADI2020}, \citet{BABAZADEH2019}, \citet{LWIN2017}, \citet{branda2016}, \citet{feng2015}, \citet{Gambrah2014}, \citet{CUI2013}, \citet{wozabal2012}, \citet{wozabal2010}, \citet{BENATI2007}, \citet{gaivoronski2005}, \citet{larsen2002}, and \citet{CAMPBELL2001} - have focused on actively managing portfolio $\text{\gls{VaR}}_{\alpha}$ rather than variance. 
This shift in risk measure can be accommodated by reformulating  the \gls{MV} model presented in \eqref{eq:markowitz} as  
\begin{equation}\label{eq:quantile_constraint}
    \begin{aligned}
    \underset{\mathbf{w} \in \Delta}{\text{maximize}} \quad & \mathbf{w}^\top \boldsymbol{\mu} \\
    \text{subject to} \quad & -\text{\gls{VaR}}_{\alpha} \left(\mathbf{w}^\top \boldsymbol{\xi} \right) \leq -r_{\text{min}}
    \end{aligned}
\end{equation}
where $\boldsymbol{\mu}:= \begin{bmatrix} \mu_1, \ldots,  \mu_n \end{bmatrix}^\top \in  \mathbb{R}^n$, and $r_{\text{min}} \in \mathbb{R}$ denotes the minimum acceptable return level, i.e., the $\text{\gls{VaR}}_\alpha$ must not fall below this threshold. 
The optimization problem can also be formulated based on a single chance constraint as
\begin{equation}\label{eq:chance_constraint}
    \begin{aligned}
    \underset{\mathbf{w} \in \Delta}{\text{maximize}} \quad & \mathbf{w}^\top \boldsymbol{\mu} \\
    \text{subject to} \quad & 1 - \mathbb{P} \left( \mathbf{w}^\top \boldsymbol{\xi} \geq r_{\text{min}} \right) \leq \alpha
    \end{aligned}
\end{equation}
where the probabilistic constraint is equivalent to the quantile-based formulation shown in \eqref{eq:quantile_constraint} [cf. \citet[p. 278]{Sarykalin2014}]. 
Independent of the representation, the altered portfolio optimization problem is difficult to solve due to the non-convexity of the risk measure and the unknown distributions of the individual investment options. 
To illustrate the former problem, the panel on the right side of Figure~\ref{fig:var_plots} presents the $\text{\gls{VaR}}_{0.05}$ as a function of the decision variable $w \in [0, 1]$ in a two-asset portfolio with finitely many scenarios. 
In this simplified example, 1,000 return observations were sampled from a multivariate normal distribution with positively correlation assets. 
The $\text{\gls{VaR}}_{0.05}$ was then computed for various weight combinations. 
In the next subsection, we discuss methodological strategies to address the complexity of chance constraints and to motivate our proposed solution approach to solve the portfolio selection problem shown in~\eqref{eq:quantile_constraint}.


\subsection{Chance Constraints: Methods and Trade-offs}\label{sec:chance_constraint_programming}


Chance-constrained programming, first introduced by \citet{Charnes1959}, is a subfield of stochastic optimization that incorporates probabilistic constraints. 
These constraints pose significant practical challenges, as the probability distributions of the underlying random variables are often unknown or difficult to evaluate [cf. \citet[p. 2222]{Pena-Ordieres2020}; \citet[p. 618]{Hong2011}]. 
Moreover, chance constraints are typically non-convex with respect to the decision variables [cf. \citet[p. 76]{Saldanha-da-Gama2024}; \citet[p. 2222]{Pena-Ordieres2020}; \citet[p. 618]{Hong2011}]. 
This further increases the difficulty of the optimization problem and makes globally optimal solutions computationally intractable without structural simplifications [cf. \citet[p. 706]{Kannan2021}]. 
To address the challenge posed by unknown probability distributions, two modeling paradigms are commonly adopted~[cf. \citet[subsection 4.4]{Saldanha-da-Gama2024}; \citet[subsection 3.2]{Birge2011}]. 
The first involves assuming parametric distributions -- typically elliptical or Gaussian -- that induce convexity in the chance constraint~[cf. \citet[p. 3]{Wang2025}; \citet[p. 2222]{Pena-Ordieres2020}]. 
This enables tractable reformulations, as demonstrated in second-order cone programs~[cf. \citet{Henrion2007}], supporting hyperplane methods~[cf. ~\citet{vanAckooij2014}], and sequential quadratic programming~[cf. \cite{Bremer2015}]. 
However, this modeling strategy will not be utilized in this paper as it is often misaligned with empirical return distributions, which exhibit skewness, heavy tails, and other stylized deviations from normality~[cf. \citet[p. 138]{zhang2018}; \citet[chapter 5]{Ruppert2011}; \citet[p. 224]{Cont2001}; \citet[pp. 394--395]{Mandelbrot1963}].


To avoid assuming specific distributional forms, the second modeling strategy relies on finite samples from the unknown distributions~[cf. \citet[p. 710]{Kannan2021}; \citet[p. 618]{Hong2011}; \citet[p. 970--971]{Nemirovski2007}]. 
These scenario sets, often derived from historical data, enable empirical approximations of the chance constraint~[cf. \citet[p. 3]{Wang2025}; \citet[p. 2222]{Pena-Ordieres2020}; \citet[p. 932]{Curtis2018}]. 
Under appropriate conditions, the resulting solutions converge almost surely to those of the true problem as the sample size increases~[cf. \citet{Pagnoncelli2009}]. 
Since this approach also aligns with established practices in financial risk management, such as historical simulation, where past return data serve as a proxy for future market behavior, it will be employed in this paper~[cf. \citet[chapter 12]{hull2023}].
Even though scenario sets simplify certain aspects of the optimization,  the empirical formulation retains the non-convexity of the original problem and may introduce artificial local minima that impair solver performance~[cf. \citet[p. 708]{Kannan2021}; \citet[p. 2237]{Pena-Ordieres2020}]. 
These effects are illustrated in the right-hand panel of Figure~\ref{fig:var_plots} and discussed in more detail by \citet{Kannan2021} and \citet{Pena-Ordieres2020}. 


Within the sample-based framework, most methods focus on computing robust local solutions with reasonable computational effort~[cf. \citet[p. 706]{Kannan2021}]. 
A common strategy is to replace the non-convex probabilistic constraint with a conservative concave surrogate~[cf. \citet{Nemirovski2007}]. 
In portfolio optimization, this is often done by substituting $\text{\gls{VaR}}_{\alpha}$ with the \gls{CVaR}~[cf. \citet{krokhmal2002}; \citet{lim2011}; \citet{norton2021}]. 
However, although \gls{CVaR} reflects expected losses in the worst-case tail~[cf. \citet[p. 1145]{Poncet2022}], it may lead to overly conservative portfolios and thus is not applied in this paper~[cf. \citet[p. 863]{wozabal2012}]. 
To avoid this problem, other sample-based approaches approximate the discontinuous indicator function that replaces the probabilistic constraint shown in~\eqref{eq:chance_constraint}~[cf. \citet[pp. 710--711]{Kannan2021}]. 
For general optimization problems involving multiple chance constraints, \citet[p. 620]{Hong2011} already proposed a \gls{DC} representation of the indicator function, while \citet[p. 721]{Kannan2021} explored smooth approximations that enable gradient-based optimization. 
A major limitation of this approach is that the probability function may exhibit flat or gully-shaped regions, which can destabilize optimization~[cf. \citet[p. 734]{vanAckooij2014b}; \citet[pp. 223--224]{Mayer2000}; \citet[pp. 2224--2225]{Pena-Ordieres2020}].
Due to this fact, in our proposed solution approach a different strategy is utilized to approximate the chance-constraint, which will be presented next.


To address the numerical difficulties of the probability function, recent work has focused on approximating the quantile-based formulation shown in~\eqref{eq:quantile_constraint}~[cf. \citet[p. 3]{Wang2025}]. 
For instance, \citet{Pena-Ordieres2020} proposed a smooth approximation of the empirical quantile function, while \citet{Wang2025} introduced a \gls{DC} formulation for joint chance constraints. 
The latter generalizes the \gls{DC} representation by \citet[p. 381]{wozabal2010}, who derived it specifically for the portfolio selection problem at hand by exploiting the relationship between $\text{\gls{VaR}}_{\alpha}$ and \gls{CVaR} under discrete distributions. 
A more in-depth discussion of this modeling framework is provided by \citet[Chapter 3]{Chandra2022}.
This decomposition enables the use of \gls{DC} programming, which is a powerful framework for addressing non-convex optimization problems~[cf. \citet{leThi2024}]. 
Subsequently, \citet{wozabal2012} applied the \gls{DCA} to this setting and proved finite convergence to a local minimum. 
Empirically, the method yields solutions close to the efficient frontier for small to medium-sized data sets~[cf. \citet[section 4]{wozabal2012}]. 
However, for more challenging $\text{\gls{VaR}}_{\alpha}$ constraints, the algorithm may fail to find feasible solutions~[cf. \citet[section 4]{wozabal2012}]. 
Moreover, it can be slow due to a high number of iterations required for larger data sets, and the resulting portfolio profits may no longer be close to the efficient frontier~[cf. \citet[section 4]{wozabal2012}]. 
Overall, these results underscore the potential of \gls{DC} programming for the optimization problem at hand, but also emphasize that the algorithmic robustness, reliability and efficiency needs to be improved, especially under demanding risk thresholds.
The next subsection briefly discusses the main idea of the \gls{DCA} and proposes our novel approach to address these research gaps.


\subsection{DC Programming Algorithms}


The \gls{DCA} framework was originally developed by \citet{Tao1988} more than 35 years ago.
The main idea of the algorithm is to iteratively solve a sequence of convex optimization problems that are easier to solve than the original problem. 
At each iteration, the non-convex objective function in \gls{DC} representation is approximated by replacing its concave part with an affine minorant derived from a subgradient [cf. \citet[p. 540]{leThi2024}]. 
This procedure results in a convex approximation of the original problem that can then be minimized efficiently based on tools from convex optimization [cf. \citet[chapter~4]{boyd2009}].  
Thus, if many iterations are required for convergence, the algorithm’s practical performance heavily depends on the scalability of the convex solver used to handle the approximate subproblems [cf. \citet[p. 53]{leThi2018}]. 


To enhance the computational efficiency of the \gls{DCA}, \citet{aragon2018} proposed the \gls{BDCA} and analyzed its convergence rate to a critical point depending on the \L ojasiewicz property. 
The key innovation of this algorithmic extension can be attributed to a line search step -- based on an Armijo-type rule with a backtracking procedure -- following each \gls{DCA} iteration. 
This enhancement enables provably larger descent steps and often results in faster convergence in practice, particularly in high-dimensional settings. 
Moreover, empirical results indicate that the line search increases the ability to escape from local minima with large suboptimality gaps [cf. \citet[p. 986]{aragon2020}]. 
While the original theoretical foundation of \gls{BDCA} was restricted to unconstrained smooth \gls{DC} problems, recent work by \citet{aragon2020} and \citet{aragon2022} has extended its applicability to certain types of non-smooth \gls{DC} problems and to \gls{DC} programs with linear inequality constraints, respectively. 
However, so far the theoretical framework of the approach has not been extended to (i) \gls{DC} programs with linear equality and inequality constraints, and (ii) \gls{DC} programs where the left part of the \gls{DC} decomposition is non-smooth. 
Moreover, the convergence rate of the algorithm depends on the specific problem structure at hand, and must be derived individually based on the \gls{KL} property.  
Overall, these findings highlight further research gaps that this paper seeks to address by custimizing the \gls{BDCA} framework to the portfolio selection problem shown in \eqref{eq:quantile_constraint}.


\subsection{Paper Contributions \& Outline}


In the line with the identified research gaps from the previous subsections, this paper proposes the \gls{BDCA} to approximately solve the optimization problem in \eqref{eq:quantile_constraint} more robustly than existing approaches.
For this purpose, we establish a novel line search scheme that enables the algorithm to be applied to \gls{DC} programs in which both components are non-smooth.
In comparison to \citet{wozabal2010}, we provide an altered \gls{DC} decomposition that ensures the strong convexity of both components. 
Based on the practical problem at hand, we then extend the existing convergence theory of the \gls{BDCA} to \gls{DC} programs with both linear equality and inequality constraints, and deduce a rigorous proof that the \gls{BDCA} with the new line search globally converges to a \gls{KKT} point with linear rate.
To obtain this result, we show that the objective corresponds to a \gls{KL} function with exponent 1/2.
In the practical part of this paper, we conduct extensive numerical experiments to assess the sensitivity of the \gls{BDCA} to its input parameters, and propose a sequential penalization of the $\text{\gls{VaR}}_{\alpha}$ constraint. 
Furthermore, we compare the proposed algorithm with the performance of \gls{DCA} and established solvers from chance-constrained programming based on four real-world data sets. 
Overall, our experiments highlight that the \gls{BDCA} outperforms the other algorithms in terms of reliability and solution quality.
In contrast to the other approaches, \gls{BDCA} demonstrates the most robust performance, as it consistently identifies feasible solutions even under the most challenging $\text{\gls{VaR}}_{X, 0.05}$ constraints and under adverse initialization.  
Ultimately, this paper provides a practical guide for transparent and easily reproducible comparisons of the \gls{VaR}-constrained portfolio selection problem using \gls{DC} programming in Python, since all code implementations are made available via GitHub.


The rest of this paper is structured as follows. 
In Section \ref{sec:var_as_dc}, the \gls{DC} representation of the discrete $\text{\gls{VaR}}_{\alpha}$ is derived based on the \gls{CVaR}. 
This is followed by the presentation of a corresponding linearly constrained \gls{DC} program. At the beginning of Section~\ref{sec:metho}, the theoretical set-ups of the \gls{DCA} and \gls{BDCA} are then introduced. 
Afterwards, the linear convergence to a \gls{KKT} point of \gls{BDCA} for the portfolio selection problem at hand is derived. 
This convergence analysis is followed by the practical application in Section~\ref{sec:prac_app}, where both \gls{DC} algorithms are compared with competing methods from chance constrained programming based on four real-world data sets. 
Ultimately, the last section summarizes all results.


\section{ Value-at-Risk as Difference of Convex Functions}\label{sec:var_as_dc}



\begin{figure}[H]%
    \centering
    \includegraphics[width=1\textwidth]{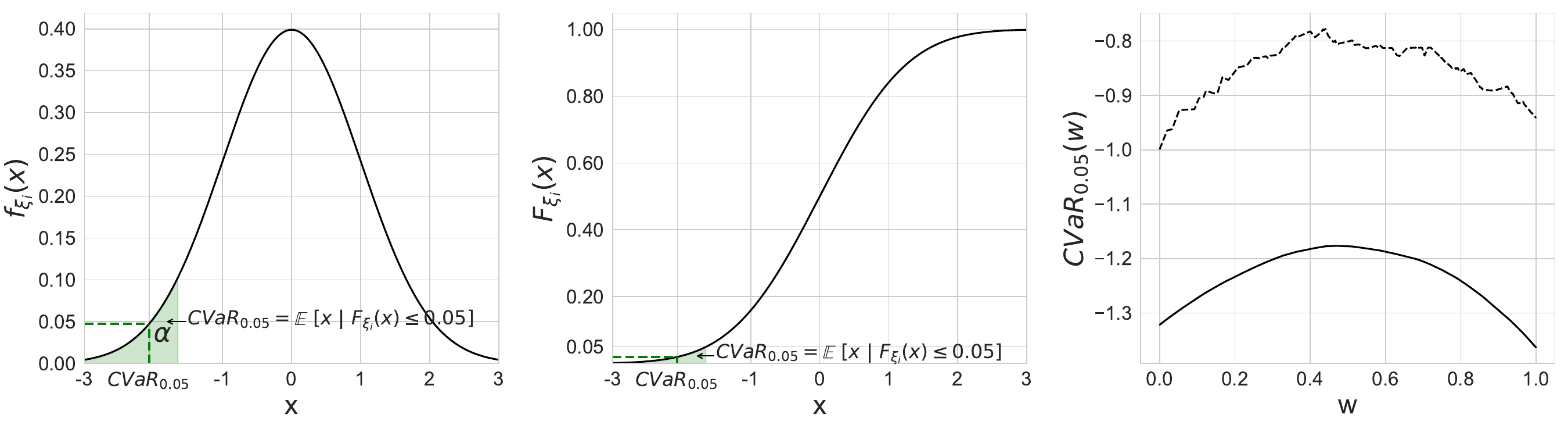}
    \captionsetup{justification=centering}
    \caption{Basic examples designed to explain the concept of $\text{\gls{CVaR}}_{\alpha}$.}%
    \label{fig:cvar_plots}%
\end{figure}

In this section, we present a natural \gls{DC} representation of $\text{\gls{VaR}}_{\alpha}$ and subsequently provide~a~corresponding linearly constrained \gls{DC} program. 
For this purpose, we extend the intuitive definition of $\text{\gls{CVaR}}_{\alpha}$ introduced earlier by providing a mathematical formulation. 
Specifically, we define $\text{\gls{CVaR}}_{\alpha}$ of a random return variable $\xi_i$ as the integral of its quantile function from 0 to $\alpha$, divided by $\alpha$. 
This enables a formal expression of the relationship between $\text{\gls{CVaR}}_{\alpha}$ and $\text{\gls{VaR}}_{\alpha}$ as 
\begin{equation}\label{eq:cvar}
    \text{\gls{CVaR}}_{\alpha} \left(\xi_i \right) := \frac{1}{\alpha} \cdot \int_{0}^{\alpha} F^{-1}_{\xi_i} \left (u \right) \ \mathrm{d}u = \frac{1}{\alpha} \cdot \int_{0}^{\alpha} \text{\gls{VaR}}_u \left(\xi_i \right) \  \mathrm{d}u,
\end{equation}
where $F^{-1}_{\xi_i}$ again denotes the inverse \gls{CDF} of $\xi_i$. 
To graphically compare the two risk measures, Figure~\ref{fig:cvar_plots} provides basic examples based on 1,000 random realizations from a standard normal distribution. 
On the left, green dashed lines indicate $\text{\gls{CVaR}}_{0.05}$ as derived from the \gls{PDF} ($f_{\xi_i}$), and in the center panel, as derived from the \gls{CDF} ($F_{\xi_i}$) of a standard normal distribution. 
The right panel displays $\text{\gls{CVaR}}_{0.05}$ as a function of the asset weight $w$ for the same portfolio example used in Figure~\ref{fig:var_plots}. 
Compared to the $\text{\gls{VaR}}_{0.05}$ function (dashed line), it is evident that $\text{\gls{CVaR}}_{\alpha}$ is concave with respect to $w$. 
Due to this concavity, $\text{\gls{CVaR}}_{\alpha}$ is a popular alternative to replace the $\text{\gls{VaR}}_{\alpha}$ constraint, as the resulting optimization problem can be efficiently solved using tools from convex optimization. 
Moreover, all panels in Figure~\ref{fig:cvar_plots} illustrate that the $\text{\gls{CVaR}}_{\alpha}$ serves as a lower bound for the $\text{\gls{VaR}}_{\alpha}$. This implies that portfolios constructed under the concave alternative will also satisfy the corresponding $\text{\gls{VaR}}_{\alpha}$ constraint.


\vspace{-0.2cm}


After outlining the relationship between $\text{\gls{CVaR}}_{\alpha}$ and $\text{\gls{VaR}}_{\alpha}$, we now explain how this result can be used to derive a \gls{DC} representation for the discrete $\text{\gls{VaR}}_{\alpha}$. 
As outlined in the introduction, it is common practice in chance-constrained programming to assume that a (representative) discrete sample can be generated either from the unknown distributions of the random variables or from historical data [cf. \citet[p. 3]{Jiang2016}]. 
For the portfolio optimization problem, this approach has the advantage that the random portfolio returns $X = \mathbf{w}^\top \boldsymbol{\xi}$ can be represented by a discrete distribution with a finite number of scenarios $S \in \mathbb{N}$. 
Specifically, the distribution of~$X$ is characterized by the set of outcomes $\{r_j\}_{j=1}^S$, where each realization is given by $r_j = \mathbf{w}^\top \mathbf{x}_j$ and occurs with probability $p_j \in [0, 1]$, such that $\sum_{j=1}^S p_j = 1$. 
This implies that the \gls{CDF} of~$X$ can be expressed as a weighted sum, and the $\text{\gls{CVaR}}_{\alpha}$ from~\eqref{eq:cvar} simplifies to
\begin{equation}\label{eq:discrete_CVaR}
    \text{\gls{CVaR}}_{X, \alpha}(\mathbf{w}) = \frac{1}{\alpha} \cdot \sum_{j=1}^{k} p_{j} \cdot \mathbf{w}^\top \mathbf{x}_{j}  + \frac{\varepsilon}{\alpha} \cdot  \mathbf{w}^\top \mathbf{x}_{k + 1},
\end{equation}
where, without loss of generality, the realizations are assumed to be sorted in increasing order, i.e., $\mathbf{w}^\top \mathbf{x}_{1}  \leq \mathbf{w}^\top \mathbf{x}_{2} \leq \cdots \leq \mathbf{w}^\top \mathbf{x}_{S}$. 
The term $\mathbf{w}^\top \mathbf{x}_{k + 1}$ is referred to as the boundary scenario, and the additional quantities $\varepsilon$ and $k$ are defined as
\begin{align}
    \begin{split}
        k &:= \underset{k}{\max} \ \sum_{j=1}^k p_{j} < \alpha \quad \text{and} \quad \varepsilon := \alpha - \sum_{j=1}^{k} p_{j},
    \end{split}
\end{align}
i.e., $k$ defines the last scenario for which the cumulative sum of the probabilities remains below~$\alpha$ and $\varepsilon$ represents the weight assigned to the boundary scenario. 
Note that for the case $p_{1} > \alpha$, the maximum of the empty set is defined as zero such that in this particular situation $\varepsilon = \alpha$ and $k = 0$ [cf. \citet[p. 865--867]{wozabal2012}].
Moreover, since $\text{\gls{CVaR}}_{X, \alpha}$ is expressed as a weighted sum of linear functions in $\mathbf{w}$ over a finite set of scenarios, it follows that the discrete $\text{\gls{CVaR}}_{X, \alpha}$ is a piecewise linear and concave function in $\mathbf{w}$ [cf. \citet[pp. 863; 870]{wozabal2012}].


The discrete $\text{\gls{CVaR}}_{X, \alpha}$ formulation in \eqref{eq:discrete_CVaR} can be used to define the discrete $\text{\gls{VaR}}_{X, \alpha}$ as a difference of two concave functions, which is given by
\begin{align}
    \text{\gls{VaR}}_{X, \alpha}(\mathbf{w})  &= \frac{\alpha}{\gamma} \cdot \text{\gls{CVaR}}_{X, \alpha}(\mathbf{w})  - \frac{\alpha - \gamma}{\gamma} \cdot \text{\gls{CVaR}}_{X, \alpha - \gamma}(\mathbf{w}) \label{eq:DC_form} \\
    &= \biggl ( \frac{1}{\gamma} \cdot \sum_{j=1}^{k} \mathbf{w}^\top \mathbf{x}_{j} \cdot p_{j} + \frac{\varepsilon}{\gamma} \cdot \mathbf{w}^\top \mathbf{x}_{k + 1}  \biggr) - \biggl( \frac{1}{\gamma} \cdot \sum_{j=1}^{k} \mathbf{w}^\top \mathbf{x}_{j} \cdot p_{j} + \frac{\varepsilon - \gamma}{\gamma} \cdot \mathbf{w}^\top \mathbf{x}_{k + 1} \biggr)  \\
    & = \frac{\varepsilon}{\gamma} \cdot \mathbf{w}^\top \mathbf{x}_{k + 1} - \frac{\varepsilon - \gamma}{\gamma} \cdot \mathbf{w}^\top \mathbf{x}_{k + 1} \\
    &= \mathbf{w}^\top \mathbf{x}_{k + 1}
\end{align}
with $0 < \gamma < \varepsilon$. This highlights that $\text{\gls{VaR}}_{X, \alpha}$ may not be uniquely defined when ties occur around the boundary scenario.  In such cases, multiple scenarios could represent the $\text{\gls{VaR}}_{X, \alpha}$, causing the discrete $\text{\gls{CVaR}}_{X, \alpha}$ to be non-differentiable at those points. 


Building on the previously derived \gls{DC} representation in \eqref{eq:DC_form}, we rewrite the generic chance-constrained maximization problem from~\eqref{eq:quantile_constraint} in empirical form with a \gls{DC} constraint. The resulting optimization problem in standard form is given by
\begin{align}\label{eq:discrete_portfolio_problem}
    \begin{split}
        \underset{\mathbf{w} \in \Delta}{\text{minimize}} \ \ & f(\mathbf{w})  \\
        \text{subject to} \  &- \text{\gls{VaR}}_{X, \alpha}(\mathbf{w}) + r_{\text{min}} \leq 0,
    \end{split}
\end{align}
where the expected return is replaced by its empirical formula $f(\mathbf{w}) :=  -\sum_{j=1}^{S} \mathbf{w}^\top \mathbf{x}_{j} \cdot p_j$, and the original maximization is equivalently expressed as a minimization problem. The feasible set $\Delta$ imposes the standard simplex constraint, but could also be modified or extended to incorporate other linear constraints. Furthermore, by substituting $\text{\gls{VaR}}_{X, \alpha}$ with the right-hand side of \eqref{eq:DC_form}, the model features a linear objective function and a non-convex \gls{DC} constraint, which could  already be addressed using iterative \gls{DC} optimization methods.


To improve tractability, especially under very strict $\text{\gls{VaR}}_{X, \alpha}$ requirements, the \gls{DC} constraint can be incorporated into the objective function via an exact penalty approach [cf. \citet[p. 870]{wozabal2012}]. Following the results of \citet{LeThi1999} on exact penalization for \gls{DC} programs, the $\text{\gls{VaR}}_{X, \alpha}$ constraint in~\eqref{eq:discrete_portfolio_problem} can be expressed in penalized form as
\begin{align}\label{eq:penalized_problem}
    \begin{split}
        \underset{\mathbf{w} \in \Delta}{\text{minimize}} \ & f(\mathbf{w})  + \tau \cdot \max \bigl[- \text{\gls{VaR}}_{X, \alpha}(\mathbf{w}) + r_{\text{min}}, 0 \bigr], 
    \end{split}
\end{align}
where $\tau > 0$ is a penalty parameter that balances return maximization and constraint satisfaction. In general, this penalized formulation can enable the algorithm to find a feasible solution even when initialized with weights $\mathbf{w}$ for which the $\text{\gls{VaR}}_{X, \alpha}$ constraint is not fulfilled. The penalty term in \eqref{eq:penalized_problem} remains a \gls{DC} and can be rewritten as
\begin{align}
    \max \big[&- \text{\gls{VaR}}_{X, \alpha}(\mathbf{w}) + r_{\text{min}}, 0 \big] 
    =  \nonumber \\ &\max \biggl[-\frac{\alpha}{\gamma} \cdot \text{\gls{CVaR}}_{X, \alpha}(\mathbf{w}) + r_{\text{min}},   - \frac{\alpha - \gamma}{\gamma} \cdot \text{\gls{CVaR}}_{X, \alpha - \gamma}(\mathbf{w}) \biggl] 
     + \frac{\alpha - \gamma}{\gamma} \cdot \text{\gls{CVaR}}_{X, \alpha - \gamma}(\mathbf{w}), \label{eq:dc_penalty_term}
\end{align}
where the \gls{DC} decomposition of $\text{\gls{VaR}}_{X, \alpha}$ is substituted, and the pointwise maximum in \eqref{eq:dc_penalty_term} defines a convex function [cf. \citet[p. 80]{boyd2009}].


Ultimately, it remains to utilize the results shown in \eqref{eq:penalized_problem} and \eqref{eq:dc_penalty_term} to reformulate the optimization problem as a linearly constrained \gls{DC} program. For this purpose, we introduce two strongly convex functions $g: \mathbb{R}^n \rightarrow \mathbb{R}$ and $h: \mathbb{R}^n \rightarrow \mathbb{R}$, which are further defined as
\begin{align}
    \begin{split}\label{eq:g_and_h}
        g(\mathbf{w}) := &   \ \psi(\mathbf{w}) + f(\mathbf{w})
         + \tau \cdot \max \biggl[ -\frac{\alpha}{\gamma} \cdot \text{\gls{CVaR}}_{X, \alpha}(\mathbf{w}) + r_{\text{min}},   - \frac{\alpha - \gamma}{\gamma} \cdot \text{\gls{CVaR}}_{X, \alpha - \gamma}(\mathbf{w})\biggl],
    \end{split}
    \\
        \begin{split}\label{eq:g_and_h2}
         h(\mathbf{w}) :=  & \  \psi(\mathbf{w}) -\tau \cdot \frac{\alpha - \gamma}{\gamma} \cdot \text{\gls{CVaR}}_{X, \alpha - \gamma}(\mathbf{w}),
    \end{split}
\end{align}
where $\psi(\mathbf{w}) := \frac{\rho}{2} \cdot \lVert \mathbf{w}\rVert^2$ with $\rho \in \mathbb{R}_{> 0}$ ensures strong convexity [cf. Definition~\ref{defi:convex_function}] of both functions, enables desirable theoretical properties in Subsection~\ref{sec:Linear_Convergence} and distinguishes our decomposition from the \gls{DC} components introduced by \citet{wozabal2010}. 
Based on these components, the optimization problem  in \eqref{eq:penalized_problem} is formulated as
\begin{align}\tag{$\mathcal{P}$}\label{eq:var_dc}
    \begin{split}
        \underset{\mathbf{w} \in \Delta}{\text{minimize}} \ &  \phi \left(\mathbf{w} \right) := g(\mathbf{w}) - h(\mathbf{w}),
    \end{split}
\end{align}
where $\phi: \mathbb{R}^n \rightarrow \mathbb{R}$ is the piecewise linear objective function, and the feasible set $\Delta$  consists of linear equality and inequality constraints.  


\vspace{0.2cm}

\begin{defi}[{Convex Function {\small [cf. \citet[chapter 3]{boyd2009}}]}]{defi:convex_function}
    A function $f: \mathbb{R}^n \rightarrow \mathbb{R}$ is said to be \textbf{convex} (with modulus $\rho = 0$) or \textbf{strongly convex} (with modulus $\rho > 0$) if
        \begin{equation*}
            f(\lambda \cdot \mathbf{w}_A + (1 - \lambda) \cdot \mathbf{w}_B) \leq \lambda \cdot f(\mathbf{w}_A) + (1 - \lambda) \cdot f(\mathbf{w}_B) - \frac{\rho}{2} \cdot \lambda \cdot (1 - \lambda) \cdot \Vert \mathbf{w}_A - \mathbf{w}_B \rVert^2
        \end{equation*}
    holds for all choices of  $\mathbf{w}_A, \mathbf{w}_B \in \mathbb{R}^n$  and $\lambda \in (0, 1)$. In contrast, a function $f$ is said to be \textbf{(strongly) concave} if $-f$ is (strongly) convex.
\end{defi}


\section{DC Programming for Value-at-Risk constrained Portfolio Optimization}\label{sec:metho}



After formulating the portfolio optimization problem in \eqref{eq:quantile_constraint} as a linearly constrained \gls{DC} program, this section customizes the theoretical setup of the \gls{BDCA} to obtain approximate solutions to \eqref{eq:var_dc}. 
In particular, a novel line search framework is proposed to address the non-smoothness of both \gls{DC} components defined in \eqref{eq:g_and_h} and \eqref{eq:g_and_h2}.
Moreover, we derive that the proposed algorithm converges linearly to a \gls{KKT} point of~\eqref{eq:var_dc}, and generalize this result to a broader class of piecewise linear \gls{DC} problems involving both linear equality and inequality constraints. 
To the best of our knowledge, these results have not been established previously. 



\subsection{The Boosted Difference-of-Convex-Functions Algorithm}


The \gls{DCA} was developed by \citet{Tao1988} more than 35 years ago, and is a powerful method to approximately solve non-convex optimization problems that can be expressed as a \gls{DC}. 
Intuitively, the key algorithmic idea is to generate a sequence of problems in which the non-convex objective function is approximated by convex functions that can be minimized efficiently [cf. \citet[p. 34--35]{LeThi2005}]. 
On a more technical level, the approach is implemented by deriving one of the infinitely many \gls{DC} decompositions of the objective function and iteratively replacing the concave part with its affine minorization,  obtained via a subgradient at the respective iterate [cf. \citet[p. 540]{leThi2024}]. 
For the optimization problem at hand, this results in a sequence of $K \in \mathbb{N}$ strongly convex subproblems of the form
\begin{equation}\tag{$\hat{\mathcal{P}}_k$}\label{eq:convex_subproblem}
    \underset{\mathbf{w} \in \Delta}{\text{minimize}} \ g(\mathbf{w}) - \Bigl(h(\mathbf{w}_k) + (\mathbf{w} - \mathbf{w}_k) ^\top \mathbf{u}_k \Bigr),
\end{equation}
where $\mathbf{u}_k \in \partial h (\mathbf{w}_k)$, with $k \in \{0, \ldots, K\}$, denotes a subgradient of $h$ at the current iterate~$\mathbf{w}_k$. 
At each iteration, these subproblems can then be solved efficiently using standard tools from convex optimization to obtain the next iterate [cf. \citet[chapter 4]{boyd2009}]. 
Due to this simple yet powerful framework, the algorithm has already been applied to a wide range of practical applications, including portfolio optimization, supply chain management, and telecommunications [cf. \citet{geremew2018}; \citet{wozabal2012}; \citet{LeThi2005}]. 
 

For the optimization problem in \eqref{eq:var_dc}, the general procedure of the \gls{DCA} is shown in Algorithm~\ref{alg:dca}. 
At each iteration, the key algorithmic steps are the (i) selection of a subgradient, (ii) minimization of the convex subproblem, and (iii) verification of the stopping condition.
The iterative procedure terminates when the unique solution of the strongly convex subproblem $\mathbf{w}_{k + 1}$ coincides with the previous iterate $\mathbf{w}_{k}$, i.e., $\lVert \mathbf{d}_k \rVert = \lVert \mathbf{w}_{k + 1} - \mathbf{w}_{k}\rVert = 0$.
This indicates that the sequence $\{\mathbf{w}_k\}$ has converged to a limit point $\bar{\mathbf{w}}$, which satisfies the \gls{KKT} conditions of problem~\eqref{eq:var_dc} [cf. \citet[p. 36]{LeThi2005}]. In other words, there exist Lagrange multipliers $\pi_1, \pi_2, \ldots, \pi_n \in \mathbb{R}_+$ and $\nu \in \mathbb{R}$ such that the following system can be constructed
\begin{align}\label{eq:KKT_conditions1}
    \begin{cases}
        \mathbf{0}_n \in \partial g(\bar{\mathbf{w}}) - \partial h(\bar{\mathbf{w}}) - \sum_{i=1}^{n}  \pi_i \cdot \mathbf{e}_i + \nu \cdot \mathbf{1}_n, \\[0.25em]
        0 = - \pi_i \cdot  \bar{\mathbf{w}}^\top  \mathbf{e}_i, \ \text{for all} \ i \in \{ 1, \ldots n\}, \\[0.25em]
        \pi_i \geq 0, \ - \bar{\mathbf{w}}^\top  \mathbf{e}_i  \leq 0, \ \text{for all} \ i \in \{ 1, \ldots, n\}, \\[0.25em]
        \mathbf{1}_n^\top \bar{\mathbf{w}}  = 1,
    \end{cases}
\end{align}
where $\mathbf{e}_i \in \{0, 1\}^n$  denotes the $i^{th}$ standard basis, and $\mathbf{1}_n := \begin{bmatrix} 1, \ldots, 1\end{bmatrix}^\top \in \{1\}^n$. 
While it is an established result that, for general \gls{DC} programs, the limit point is approached asymptotically, the \gls{DCA} exhibits finite convergence when applied to polyhedral \gls{DC} programs, as proposed by \citet{wozabal2012}  for the problem at hand [cf. \citet[p. 36]{LeThi2005}].
In theory, the latter setup guarantees that $\lVert \mathbf{d}_k \rVert$ will eventually be exactly zero. However, in practice, the number of iterations required to reach this result is not known a priori and may grow substantially with the number of piecewise-linear terms. 
For instance, \citet[p. 879]{wozabal2012} observed that the polyhedral \gls{DCA} becomes computationally impractical for large-scale portfolio examples, as more iterations are required for convergence, making its performance heavily dependent on the scalability of the underlying convex solver [cf. \citet[p. 53]{leThi2018}].


\begin{algorithm}[H]
\DontPrintSemicolon
  \KwInput{Initial $\mathbf{w}_0 \in \Delta, K \in \mathbb{N}$, $\tau \in \mathbb{R}_{>0}$ and $\rho \in \mathbb{R}_{>0}$}
  $k \gets 0$ \\
  \For{$k \leq K$}
   {
    \vspace{0.1cm}
    $\mathbf{u}_k \gets$ random element of $\partial h(\mathbf{w}_k)$\\
    \vspace{0.1cm}
    $\mathbf{w}_{k + 1} \gets \underset{\mathbf{w} \in \Delta}{\text{arg min}} \ \ g(\mathbf{w}) - \mathbf{w}^\top \mathbf{u}_k$ \\
    \vspace{0.1cm}
    $\mathbf{d}_k \gets \mathbf{w}_{k + 1} - \mathbf{w}_{k}$ \\
    \vspace{0.15cm}
    \If{$\lVert \mathbf{d}_k \rVert = 0$}
    {
    \vspace{0.1cm}
    \textbf{break}
    }
    \vspace{0.1cm}
    $k \gets k + 1$ \\
   }
   \vspace{0.1cm}
   \KwOutput{$\mathbf{w}_k$}
   \vspace{0.1cm}
\caption{\acrfull{DCA}}\label{alg:dca}
\end{algorithm}


To enhance the computational efficiency of the \gls{DCA}, \citet{aragon2018} proposed the \gls{BDCA}. 
The key innovation of this algorithmic extension is the incorporation of a line search step following each \gls{DCA} iteration.
Theoretically, this procedure results in provably larger descent steps and often faster convergence in practice [cf. \citet[Corollary 1] {aragon2018}].  
Moreover, empirical results indicate that the line search increases the ability to escape from poor local minima, i.e., points far from the global solution [cf. \citet[p. 986]{aragon2020}]. 
While the theoretical foundation of the original \gls{BDCA} was limited to unconstrained smooth \gls{DC} problems, \citet{aragon2020} and \citet{aragon2022} extended its applicability to \gls{DC} problems with non-smooth $h$ and to inequality-constrained programs, respectively. 
However, since the optimization problem shown in~\eqref{eq:var_dc} corresponds to a \gls{DC} program (i) with linear equality and inequality constraints, and (ii) where both components of the \gls{DC} decomposition are non-smooth, the previous extensions are not sufficient to apply the \gls{BDCA}. 
The current setup of the line search procedure is particularly problematic as it requires the function $g$ to be smooth in order to guarantee a finite termination.  
To overcome this limitation, we propose an alternative line search scheme which, to the best of our knowledge, has not previously been considered in the context of the \gls{BDCA}.
The resulting algorithmic procedure, customized for problem \eqref{eq:var_dc}, is presented in Algorithm~\ref{alg:bdca}, where the key components of the line search are the
\begin{enumerate}
    \item computation of a search direction $\mathbf{d}_k: = \mathbf{y}_k - \mathbf{w}_k$,
    \item determination of the step size $\lambda_k \in \mathbb{R}_+$ based on the boundary of $\Delta$, and
    \item application of an Armijo-type condition combined with a backtracking procedure.
\end{enumerate}
At each iteration, the line search is initialized on the boundary of $\Delta$ to simplify the selection of~$\lambda_k$. 
If the $\text{\gls{VaR}}_{X, \alpha}$ constraint is not violated throughout the search direction~$\mathbf{d}_k$, the objective function behaves linearly along this path. 
As a result, starting at the boundary leads to the largest possible reduction in the objective when $\mathbf{d}_k$ is a descent direction. 
Computationally, this can reduce runtime, as larger steps may be accepted without repeatedly scaling down $\lambda_k$. 
In cases involving challenging $\text{\gls{VaR}}_{X, \alpha}$ constraints, this strategy combined with a large $\beta$ can also help escape infeasible regions or identify ``islands'' in the potentially disconnected feasible domain, which will be demonstrated later in this subsection. 
For the problem at hand, starting from the boundary can also encourage sparser solutions that allocate the portfolio weights across fewer assets.

\vspace{0.2cm}


\begin{algorithm}[H]
\DontPrintSemicolon
  \vspace{0.1cm}
  \KwInput{Initial $\mathbf{w}_0 \in \Delta$, $K \in \mathbb{N}$, $\tau \in \mathbb{R}_{>0}$, $\rho \in \mathbb{R}_{>0}$, and $\beta \in (0, 1)$}
  \vspace{0.1cm}
  $k \gets 0$ \\
  \vspace{0.1cm}
  \For{$k \leq K$}
   {
   \vspace{0.1cm}
   $\mathbf{u}_k \gets$ random element of $\partial h(\mathbf{w}_k)$\\
   \vspace{0.1cm}
    $\mathbf{y}_k \gets \underset{\mathbf{w} \in \Delta}{\text{arg min}} \ \ g(\mathbf{w}) - \mathbf{w}^\top \mathbf{u}_k$ \\
    \vspace{0.1cm}
    $\mathbf{d}_k \gets \mathbf{y}_k - \mathbf{w}_k$ \\
    \vspace{0.1cm}
    \If{$\lVert \mathbf{d}_k \rVert =0$}
    {\vspace{0.1cm}
    \textbf{break}
    }
    \Else
    {
        \vspace{0.1cm}
        $\lambda_k \gets  \underset{( \mathbf{w}_k + \lambda \cdot \mathbf{d}_k  ) \in \Delta}{\text{arg max}} \  \lambda$ \label{line:lam_boundary_set}\\
        \vspace{0.1cm}
        \While{$\phi(\mathbf{w}_k + \lambda_k \cdot \mathbf{d}_k) > \phi(\mathbf{w}_k) - \rho \cdot  \lambda_k^2 \cdot \lVert \mathbf{d}_k \rVert^2$ \label{line:armijo_condition}}
        {
            \vspace{0.1cm}
            $\lambda_k \gets \max \bigl \{ 1, \lambda_k \cdot \beta \bigr \}$ \label{line:backtracking}
        }
     }
 
     \vspace{0.1cm}
    $\mathbf{w}_{k + 1} \gets \mathbf {w}_k + \lambda_k \cdot \mathbf{d}_k$ \\
    \vspace{0.1cm}
    }
   \KwOutput{$\mathbf{w}_k$}
   \vspace{0.1cm}
\caption{\acrfull{BDCA}}\label{alg:bdca}
\end{algorithm}


\vspace{0.2cm}

\begin{rema}[Line Search Procedure]{rema:new_line_search}
    The line search shown in Lines \ref{line:lam_boundary_set} to \ref{line:backtracking} of Algorithm~\ref{alg:bdca} is 
    similar to that \citet{FUKUSHIMA1981}.
    However, in their scheme the function $h$ must be differentiable and the update of the solution vector is performed based on $\mathbf{w}_{k+1} = \mathbf{w}_{k} + \lambda_{\text{FM}}^l \cdot \mathbf{d}_k$ where the step size is limited to $\lambda_{\text{FM}} \in (0, 1)$ and $l \in \mathbb{N}$.
    Thus, $\lambda_{\text{FM}}^l$ is always accepted smaller than $1$, while our step size $\lambda_k \ge 1$ for all $k$, and the procedure leads to totally different iterations. 
    
    \noindent In comparison to \citet{aragon2022}, the Armijo-type condition has changed as in their line search the condition $\phi \left( \mathbf{y}_k + \lambda_k \cdot \mathbf{d}_k \right) = \phi \left (\mathbf{y}_k \right ) - \alpha \cdot \lambda_k^2 \cdot \lVert \mathbf{d}_k \rVert^2$ with $\alpha > 0$ must be fulfilled and the update is performed based on $\mathbf{w}_{k+1} = \mathbf{y}_{k} + \lambda_k \cdot \mathbf{d}_k$.
    Moreover, it must be verified that $g$ is differentiable at $\mathbf{y}_k$, as otherwise $\mathbf{d}_k$ might be an ascent direction and the line search could become an infinite loop [cf. \citet[Remark~1]{aragon2018}].
    This limitation no longer applies to the new procedure presented in Algorithm~\ref{alg:bdca}, as the line search terminates in the worst case when $\lambda_k = 1$, yielding the solution of the convex subproblem: $\mathbf{w}_{k+1} = \mathbf{w}_k + \mathbf{d}_k = \mathbf{w}_k + \left(\mathbf{y}_k - \mathbf{w}_k \right)$. 
    If the step size is accepted with $\lambda_k > 1$, then the updates of the two line search procedures remain similar.

\end{rema}



\begin{figure}[H]%
    \centering
    \includegraphics[width=1\textwidth]{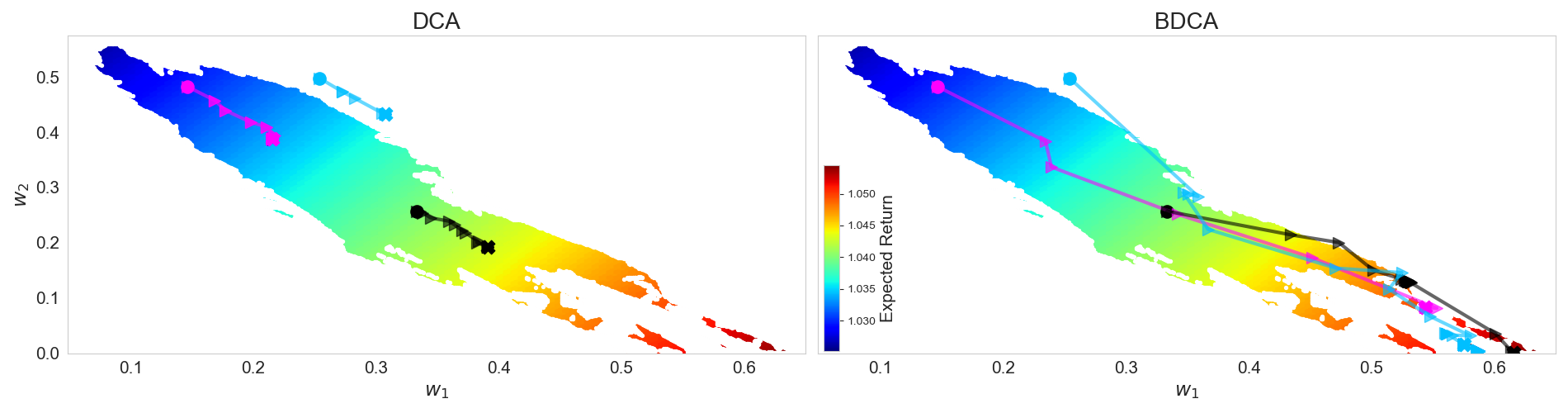}
    \caption{Optimization trajectory comparison between \gls{DCA} and \gls{BDCA} based on synthetic data example. 
    The horizontal and vertical axes represent the weights $w_1$ and $w_2$, respectively, while the color shading indicates the expected portfolio profit for combinations that satisfy the $\text{\gls{VaR}}_{X, 0.05}$ and $\Delta$ constraints.}%
    \label{fig:heat_map_dca_bdca}%
\end{figure}

After introducing the \gls{BDCA} with the new line search framework, we demonstrate its benefits over the \gls{DCA} based on a first toy example.  
For this purpose, we generate $2{,}000$ scenarios using a multivariate normal distribution with covariance matrix
\begin{equation*}
\boldsymbol{\Sigma} = \begin{bmatrix}
                        1.0 & 0.7 & 0.3 \\
                        0.7 & 1.0 & 0.4 \\
                        0.3 & 0.4 & 1.0
                        \end{bmatrix},
\end{equation*}
and mean vector $\boldsymbol{\mu} = \begin{bmatrix} 1.05 & 1.03 & 1.055 \end{bmatrix}^\top$.
Given this random sample, Figure~\ref{fig:heat_map_dca_bdca} visualizes the optimization trajectories of both algorithms, initialized from three randomly chosen starting points. 
The horizontal and vertical axes of Figure~\ref{fig:heat_map_dca_bdca} represent different combinations of the weights $w_1$ and $w_2$, respectively, while the value of the third weight can be derived as $w_3 = 1 - w_1 - w_2$.
The color shading highlights the expected portfolio profit for weight combinations that satisfy both the $\text{\gls{VaR}}_{X, 0.05}$ and $\Delta$ constraints, with infeasible regions shown in white.
The right subplot illustrates how \gls{BDCA} effectively escapes low-return regions and quickly advances toward more favorable areas. 
Moreover, the algorithm is able to detect feasible islands located in the bottom right of the non-convex solution space. 
In contrast, the left subplot shows that the \gls{DCA} frequently becomes trapped in suboptimal regions and progresses more slowly. 
This is also visually apparent in the optimization trajectories, where the dots representing successive iterates of \gls{DCA} are densely clustered, while the steps of \gls{BDCA} are generally more widely spaced.
Overall, these first empirical findings underscore the usefulness of the proposed line search step. 


\vspace{0.1cm}

\vspace{-0.2cm}


In line with the previous empirical findings, the advantages of the line search can also be substantiated from a theoretical point of view.
To this end, Theorem~\ref{theo:main1} summaries key properties of Algorithm~\ref{alg:bdca}.
The first result establishes the improvement achieved by solving the strongly convex subproblem, while the second confirms that an even longer step size is used in all iterations where the line search procedure terminates with $\lambda_k > 1$.
These results then form the basis for deriving that the sequence of objective function values is monotonously decreasing and convergent to a limit point, as summarized in the third bullet point. 
Ultimately, the fourth result also deduces that the algorithm converges to a \gls{KKT} point of \eqref{eq:var_dc}, which serves as the theoretical basis for deriving the convergence rate in the next subsection.


\vspace{0.2cm}

\begin{theo}{theo:main1}
The sequences $\{ \mathbf{w}_k \}$ and $\{ \mathbf{y}_k \}$ generated by Algorithm \ref{alg:bdca} satisfy the following properties 
    \begin{itemize}
    \item [(i)] For all $k \ge 0$, the strongly convex subproblem yields
    $$\phi(\mathbf{y}_k) \le \phi(\mathbf{w}_k) - \rho \cdot \lVert \mathbf{d}_k\rVert^2.$$
    \item [(ii)] The line search in Lines \ref{line:lam_boundary_set} to \ref{line:backtracking} terminates in finite steps, and guarantees
    $$\phi(\mathbf{w}_{k+1}) 
         \le \phi(\mathbf{w}_k) - \rho \cdot \lambda_k^2\cdot \lVert \mathbf{d}_k\rVert^2.$$
    \item[(iii)] The sequence $\bigl \{\phi \left( \mathbf{w}_k \right) \bigr\}$ is monotonously decreasing and convergent to some $\Bar{\phi}$.
    \item [(iv)] There exists a sub-sequence of $\bigl \{ \mathbf{w}_k \bigr\}$ converging to a \gls{KKT} point of \eqref{eq:var_dc}.
    \end{itemize}
\end{theo}


\begin{prf}{prf:XX}

    \begin{itemize}
        \item[(i)] 
        By the first order optimality condition of the convex subproblem \eqref{eq:convex_subproblem}, we have
        $$
        \mathbf{0}_n \in \partial g(\mathbf{y}_k)- \mathbf{u}_k + N_{\Delta}(\mathbf{y}_k) \quad \text{for all} \quad k \geq 0,
        $$
        where $N_{\Delta}(\mathbf{y}_k) := \bigl\{ \mathbf{v} \in \mathbb{R}^n: \mathbf{v}^\top(\mathbf{y}-\mathbf{y}_k) \le 0 \quad \text{for all } \mathbf{y} \in \Delta \bigr\}$ is the normal cone of $\Delta$ at $\mathbf{y}_k$.
        This implies that there exists $\mathbf{v}_k \in N_{\Delta}(\mathbf{y}_k)$ such that
        $$
        \mathbf{u}_k -\mathbf{v}_k \in \partial g \left(\mathbf{y}_k \right),
        $$
        which can be combined with the strong convexity of $g$ to obtain
        \begin{align} \label{g_convexity}
        g(\mathbf{w}_k) - g(\mathbf{y}_k) 
        &\ge  (\mathbf{u}_k -\mathbf{v}_k)^\top (\mathbf{w}_k - \mathbf{y}_k) +\frac{\rho}{2} \cdot \lVert\mathbf{y}_k - \mathbf{w}_k\|^2 \nonumber \\
        & = \mathbf{u}_k^\top (\mathbf{w}_k - \mathbf{y}_k) 
        -\mathbf{v}_k^\top (\mathbf{w}_k - \mathbf{y}_k) 
        +\frac{\rho}{2} \cdot \|\mathbf{y}_k - \mathbf{w}_k\|^2 \nonumber\\
        &\ge \mathbf{u}_k^\top (\mathbf{w}_k - \mathbf{y}_k) 
        +\frac{\rho}{2} \cdot \|\mathbf{y}_k - \mathbf{w}_k\|^2,
        \end{align}
        where the last inequality follows by $\mathbf{v}_k \in N_{\Delta}(\mathbf{y}_k)$ and $\mathbf{w}_k \in \Delta$.
        On the other hand, since $\mathbf{u}_k \in \partial h(\mathbf{w}_k)$ and $h$ is $\rho$-strongly convex, we have
        \begin{align} \label{h_convexity}
        h(\mathbf{y}_k) - h(\mathbf{w}_k) \ge  \mathbf{u}_k^\top (\mathbf{y}_k - \mathbf{w}_k) +\frac{\rho}{2} \cdot \|\mathbf{y}_k - \mathbf{w}_k\|^2.
         \end{align}
        Adding \eqref{g_convexity} and \eqref{h_convexity}, we deduce
        \begin{equation}\label{eq:first_result}
        \phi(\mathbf{y}_k) \le \phi(\mathbf{w}_k) - \rho\cdot \lVert \mathbf{d}_k\rVert^2.
        \end{equation}
        
    
        \item[(ii)] 
        Since $\mathbf{y}_k \in \Delta$, the solution $\lambda_k$ of the maximization problem in Line~\ref{line:lam_boundary_set} always exists and is bounded from below by $\lambda_k \ge 1$. 
        In addition, the result in \eqref{eq:first_result} implies that the line search procedure in Lines \ref{line:armijo_condition} and \ref{line:backtracking} is well defined because in the worst case, the step size $\lambda_k = 1$ will be accepted.
        These results then confirm
        $$\phi(\mathbf{w}_{k+1}) 
         \le \phi(\mathbf{w}_k) - \rho \cdot \lambda_k^2\cdot \lVert \mathbf{d}_k\rVert^2.$$

    
        \item[(iii)]
         From the line search procedure and $\lambda_k \ge 1$, we have
         \begin{align} \label{decreasing objective}
             \phi(\mathbf{w}_{k+1}) 
         \le \phi(\mathbf{w}_k) - \rho \cdot \lambda_k^2\cdot \lVert \mathbf{d}_k\rVert^2
         \le \phi(\mathbf{w}_k) - \rho\cdot \lVert \mathbf{d}_k\rVert^2, 
         \end{align}
        which means that the sequence $\bigl \{\phi \left( \mathbf{w}_k \right) \bigr\}$ is monotonously decreasing. 
        It remains to show that the objective function is bounded from below. 
        This follows directly from the individual components in \eqref{eq:penalized_problem} that are themselves bounded from below. 
        Specifically, the penalty term $\tau \cdot \max \bigl[- \text{\gls{VaR}}_{X, \alpha}(\mathbf{w}) + r_{\text{min}}, 0 \bigr] \geq 0$ is always nonnegative. 
        In contrast, $f$ represents the negative portfolio return, which could be arbitrarily small, at least in theory.  
        However, since $f$ is continuous and defined over a compact set on $\mathbb{R}^n$, i.e., $f(\Delta)$, the Extreme-Value Theorem implies that the function attains its maximum and minimum within $\Delta$ and thus must be bounded as well.

        \item[(iv)]
        We notice that $\mathbf{w}_k \in \Delta $ for all $k\ge 0$. Since $\Delta$ is a compact set (i.e., closed and bounded), the sequence $\bigl \{ \mathbf{w}_k \bigr\}$ is bounded. As a consequence, there exists a sub-sequence $\bigl \{ \mathbf{w}_{k_t} \bigr\}$ of $\bigl  \{ \mathbf{w}_k \bigr\}$ converging to some $\bar{\mathbf{w}}$. From \eqref{decreasing objective} we deduce that
        $$
        \rho\cdot \lVert \mathbf{y}_k - \mathbf{w}_k \rVert^2 = \rho\cdot \lVert \mathbf{d}_k\rVert^2 \le \phi \left( \mathbf{w}_k \right) - \phi \left( \mathbf{w}_{k+1} \right) \to 0,
        $$
        which also implies that 
        $
        \lim_{t \to \infty} \mathbf{y}_{k_t} = \bar{\mathbf{w}}.
        $
        Note that the \gls{KKT} conditions of the subproblem \eqref{eq:convex_subproblem} is
            \begin{align*}
                \begin{cases}
                  \mathbf{0}_n \in  \partial g(\mathbf{y}_k) -  \mathbf{u}_k - \sum_{i=1}^{n} \pi_{k,i} \cdot \mathbf{e}_i +  \nu_{k} \cdot \mathbf{1}_n, \\[0.25em]
                    - \pi_{k,i} \cdot \mathbf{e}_i^\top \mathbf{y}_k = 0, \pi_{k,i} \geq 0, \ - \mathbf{e}_i^\top \mathbf{y}_k \leq 0, \ \text{for all} \ i \in \{ 1, \ldots n\}, \\[0.25em]
                    \mathbf{1}_n^\top \mathbf{y}_k = 1.
                \end{cases}
            \end{align*}
        Since $\mathbf{u}_k  \in  \partial h(\mathbf{w}_k)$, we can write the last system as 
                \begin{align}\label{eq:KKT_cond_bis}
                \begin{cases}
                  \mathbf{u}_k  \in  \partial h(\mathbf{w}_k) \cap \Bigl( \partial g(\mathbf{y}_k) -    \sum_{i=1}^{n} \pi_{k,i} \cdot \mathbf{e}_i +  \nu_{k} \cdot \mathbf{1}_n \Bigr), \\[0.25em]
                    -\pi_{k,i} \cdot \mathbf{e}_i^\top \mathbf{y}_k = 0, \pi_{k,i} \geq 0, \ - \mathbf{e}_i^\top \mathbf{y}_k \leq 0, \ \text{for all} \ i \in \{ 1, \ldots n\}, \\[0.25em]
                    \mathbf{1}_n^\top \mathbf{y}_k = 1.
                \end{cases}
            \end{align}
        Furthermore, since $h$ is continuous and $\mathbf{u}_k  \in  \partial h(\mathbf{w}_k)$ with $\mathbf{w}_k \in \Delta \subset \textbf{intdom } h$, the sequence $\{ \mathbf{u}_k\}$ is bounded. 
        Hence, we can extract a convergent subsequence and thus assume $\mathbf{u}_{k_t} \rightarrow \Bar{\mathbf{u}}$. 
        In combination with \eqref{eq:KKT_cond_bis}, we obtain the following \gls{KKT} system
            \begin{align}\label{eq:KKT_Cond2bis}
            \begin{cases}
               \mathbf{u}_{k_t} \in \partial h(\mathbf{w}_{k_t}) \cap \Bigl( \partial g(\mathbf{y}_{k_t}) - \sum_{i=1}^{n} \pi_{k_t,i} \cdot \mathbf{e}_i + \nu_{k_t} \cdot \mathbf{1}_n \Bigr) \\
                -\pi_{k_t,i} \cdot  \mathbf{e}_i^\top \mathbf{y}_{k_t} = 0, \ \pi_{k_t,i} \geq 0, \ - \mathbf{e}_i^\top \mathbf{y}_{k_t} \leq 0, \  \text{for all} \ i \in \{ 1, \ldots n\}, \\
                \mathbf{1}_n^\top \mathbf{y}_{k_t} = 1,
            \end{cases}
            \end{align}
            for which we aim to take the limit as $t \rightarrow \infty$. 
            For this purpose, it only remains to show that the subsequences of Lagrange multiplier $\{ \pi_{k_t, i}\}$ and $\{ \nu_{k_t}\}$ are also bounded. 
            This can be proven based on the dual form of the \gls{MFCQ} that can be derived based on Motzkin's Transposition Theorem [cf. \citet[p. 4]{Solodov2011}; \citet[p. 44]{MANGASARIAN1967}]. 
            For the sake of brevity, we only show that the problem in \eqref{eq:var_dc} fulfills the \gls{MFCQ}. 
            As the feasible set has only one equality constraint, i.e., $\Tilde{h}(\mathbf{w}) := \mathbf{w}^\top \mathbf{1}_n - 1 = 0$, it is trivial to show that
            \begin{equation}\label{eq:mfcq_1}
                a \cdot \nabla \Tilde{h}(\mathbf{z}) = a \cdot\mathbf{1}_n = \mathbf{0}_n  \quad \text{with} \quad a \in \mathbb{R},
            \end{equation}
            has only one solution at $a = 0$.
            Thus, the independence requirement with respect to the gradients of the equality constraints is always fulfilled. 
            By selecting the elements of $\mathbf{a}:= \begin{bmatrix} a_1, \ldots, a_n \end{bmatrix}^\top \in \mathbb{R}^n$ such that $\sum_{i=1}^n a_i = 0$ and
            \begin{align*}
                    a_i \in
                    \begin{cases}
                         \mathbb{R}_{>0} \quad &\text{if} \quad w_i = 0, \\
                        \mathbb{R} \quad &\text{if} \quad w_i > 0, \\
                    \end{cases}
            \end{align*}
            is fulfilled, we automatically can obtain the following results 
            \begin{equation}\label{eq:mfcq_2}
                 - \mathbf{e}_i^\top  \mathbf{a}  = - a_i < 0 \quad \forall i \in I(\mathbf{w}) \quad \text{and} \quad \nabla \Tilde{h}(\mathbf{w})^\top \mathbf{a} = \sum_{i=1}^n a_i = 0,
            \end{equation}
            where $I(\mathbf{w})$ denotes the set of active inequality constraints.
            The combination of \eqref{eq:mfcq_1} and \eqref{eq:mfcq_2} shows that all feasible points fulfill the \gls{MFCQ}, and thus the sequences of Lagrange multipliers are bounded, i.e.,
            \begin{equation}
                \underset{t \rightarrow \infty}{\lim} \ \pi_{k_t, i} = \bar{\pi}_i, \ \text{for all} \ i \in \{1, \ldots, n\} \quad \text{and} \quad \underset{t \rightarrow \infty}{\lim} \ \nu_{k_t} = \bar{\nu}.
            \end{equation}
            Due to the closedness of the graph of $\partial h$ and $\partial g$, we can pass to the limit $t \rightarrow \infty$ in~\eqref{eq:KKT_Cond2bis} to obtain
            \begin{align}\label{eq:KKT_Cond3}
                \begin{cases}
                  \Bar{\mathbf{u}} \in \partial h(\Bar{\mathbf{w}}) \cap \Bigl( \partial g(\Bar{\mathbf{w}})  - \sum_{i=1}^{n} \bar{\pi}_{i} \cdot \mathbf{e}_i + \bar{\nu} \cdot \mathbf{1}_n \Bigr), \\
                    - \bar{\pi}_{i} \cdot \mathbf{e}_i^\top \Bar{\mathbf{w}} = 0, \ \bar{\pi}_{i} \geq 0, \ -\mathbf{e}_i^\top, \Bar{\mathbf{w}} \leq 0, \  \text{for all} \ i \in \{ 1, \ldots n\}, \\
                    \mathbf{1}_n^\top \Bar{\mathbf{w}} = 1,
                \end{cases}
            \end{align}
            which confirms that $\Bar{\mathbf{w}}$ is a \gls{KKT} point of \eqref{eq:var_dc}.
        
    \end{itemize}

\end{prf}



\subsection{Linear Convergence for Piecewise Linear Objective Functions}\label{sec:Linear_Convergence}


After proposing a new version of the \gls{BDCA} and proving its convergence to a \gls{KKT} point, this subsection establishes the corresponding convergence rate for the problem at hand.
To establish our result, we prove that the objective function satisfies the requirements of a piecewise linear function as defined below in Definition~\ref{def:piecewise_linear}.


\vspace{0.2cm}

\begin{defi}[Piecewise Linear Function]{def:piecewise_linear}
    A function $f: \mathbb{R}^n \rightarrow \mathbb{R}$ is called \textbf{piecewise linear} if it can be expressed as
    \begin{equation*}
        f(\mathbf{w}) = \max_{i \in \{1, \ldots, m\}} \Bigl (\mathbf{c}_i^\top \mathbf{w} + b_i \Bigr )
    \end{equation*}
    for some $m \in \mathbb{N}, \mathbf{w} \in \mathbb{R}^n, \mathbf{c}_i \in \mathbb{R}^n$ and $b_i \in \mathbb{R}$.
\end{defi}

\vspace{-0.2cm}


After showing that our objective conforms to Definition~\ref{def:piecewise_linear}, we deduce based on Definition~\ref{defi:KL_fct_exponent} and Theorem~\ref{theo:kl} that piecewise linear functions are \gls{KL} functions with exponent 1/2.
To the best of our knowledge, this result has not been shown previously.
The \gls{KL} exponent plays a central role in the analysis of convergence rates for various algorithms. 
In particular, within \gls{DC} programming, it is frequently used to establish a rigorous convergence theory. 
In this paper, it enables us to deduce the linear convergence rate of Algorithm~\ref{alg:bdca}, as stated in Theorem~\ref{theo:linear_rate}.


\vspace{0.2cm}

\begin{defi}[\gls{KL} Function and Exponent]{defi:KL_fct_exponent}
    A proper closed function $f: \mathbb{R}^n \rightarrow \mathbb{R}$ is called a \textbf{\gls{KL}} \textbf{function} if for all $\Bar{\mathbf{w}} \in \text{dom} \left(\partial f \right)$ there exist a neighborhood $\mathcal{N}$ around $\Bar{\mathbf{w}}$, a constant $\eta \in (0, \infty]$, and a continuous concave function $\varphi: [0, \eta) \rightarrow \mathbb{R}_+$ with $\varphi(0) = 0$ such that:
    \begin{enumerate}
        \item $\varphi$ is continuously differentiable on $(0, \eta)$ with $\nabla \varphi(t) > 0$ for all $t \in (0, \eta)$;
        \item for all $\mathbf{w} \in \mathcal{N}$ with $f(\Bar{\mathbf{w}}) < f(\mathbf{w}) < f(\Bar{\mathbf{w}}) + \eta$, it holds that
        \[
        \nabla \varphi(f(\mathbf{w}) - f(\Bar{\mathbf{w}})) \cdot \text{dist}(0, \partial f(\mathbf{w})) \geq 1.
        \]
    \end{enumerate}
    If for all $\Bar{\mathbf{w}} \in \text{dom} \left(\partial f \right)$, $\varphi$ can be chosen as $\varphi \left(t \right) = M \cdot t^{1 - \theta_f}$ for some $M > 0$, $t >0 $ and $\theta_f \in [0, 1)$, i.e., there exist a constant $\epsilon > 0$ such that
    \[
    \text{dist}(0, \partial f(\mathbf{w})) \geq t \cdot \bigl(f(\mathbf{w}) - f(\Bar{\mathbf{w}}) \bigr)^{\theta_f},
    \]
    whenever $\| \mathbf{w} - \Bar{\mathbf{w}} \| \leq \epsilon$ and $f(\Bar{\mathbf{w}}) < f(\mathbf{w}) < f(\Bar{\mathbf{w}}) + \eta$, then $f$ is said to be a \gls{KL} function with \gls{KL} \textbf{exponent} $\theta_f$.
\end{defi}


\begin{theo}[{\citet[p. 1221]{li2018}}]{theo:kl}
    Suppose that $f: \mathbb{R}^n \rightarrow \mathbb{R}$ is a proper closed function of the form
    \begin{equation*}
        f(\mathbf{w}) = l(\mathbf{A} \mathbf{w}) + \underset{1 \leq i \leq r}{\min} P_i (\mathbf{w}),
    \end{equation*}
    where $\mathbf{A} \in \mathbb{R}^{m \times n}$, $P_i$ are proper closed polyhedral functions for $i = 1, \ldots, r$ and $l$ is a proper closed convex function with an open domain, is strongly convex on any compact convex subset of dom $l$ and is twice continuously differentiable on dom $l$. Suppose in addition that $f$ is continuous on dom $\partial f$. Then $f$ is a \gls{KL} function with an exponent of~$\frac{1}{2}$. 
\end{theo}


\begin{theo}[Linear Convergence to \gls{KKT} Point]{theo:linear_rate}
Let $\bar{\mathbf{w}} \in \Delta$  be a limit point of $\{ \mathbf{w}_k \}$. 
If $ \nabla g $ or $ \nabla h$ is locally Lipschitz continuous around $\bar{\mathbf{w}}$, then the sequence $\{\mathbf{w}_k\}$ converges linearly to $\bar{\mathbf{w}}$.
\end{theo}


\begin{prf}{prf:convergence_result}
   We first demonstrate that the objective function in~\eqref{eq:var_dc} is a \gls{KL} function with exponent~$\frac{1}{2}$. 
   For this purpose, we show that the $\text{\gls{CVaR}}_{X, \alpha}(\mathbf{w})$, as defined in \eqref{eq:discrete_CVaR}, can be rewritten in the general form of a piecewise linear function, as specified in Definition~\ref{def:piecewise_linear}. 
   To obtain this representation, we use the results from \cite{rockafellar2000}, yielding 
    \begin{align}
        \text{\gls{CVaR}}_{X, \alpha} \left( \mathbf{w} \right)
        &= \underset{a \in \mathbb{R}}{\max} \ \biggl\{ a - \frac{1}{\alpha} \cdot \mathbb{E} \Bigl( \min \left[X - a, 0 \right] \Bigr) \biggr\} \\
        &= \underset{a \in \mathbb{R}}{\max} \ \biggl\{ a - \frac{1}{\alpha} \cdot \sum_{i=1}^S \min \left[p_i \cdot \mathbf{w}^\top \mathbf{x}_i - a, 0 \right]  \biggr\} \\
        &= \underset{j \in \{1, \ldots, S\} }{\max} \ \biggl\{ \underbrace{\biggl (1 + \frac{m}{\alpha} \biggr) \cdot a_j }_{b_j}  + \mathbf{w}^\top \underbrace{\biggl (- \frac{1}{\alpha} \cdot \sum_{i \in \mathcal{Z}_{(\mathbf{w}, a_j)}} p_i \cdot \mathbf{x}_i\biggr)}_{\mathbf{a}_j}   \biggr\}
    \end{align}
    where $\mathcal{Z}_{(\mathbf{w}, a_j)} := \left \{i: \mathbf{w}_i^\top \mathbf{x}_i < a_j  \right \}$ and the optimal value of $a_j \in \mathbb{R}$ at the maximum corresponds to the discrete $\text{\gls{VaR}}_{X, \alpha}(\mathbf{w})$. 
    Hence, the maximum is attained over only $S$ candidate values, and $\text{\gls{CVaR}}_{X, \alpha}(\mathbf{w})$ can be expressed as the maximum over $S$ affine functions in $\mathbf{w}$. 
    
    As second step, we prove that $\phi$ fulfills Definition~\ref{def:piecewise_linear} as it consists of the sum of a linear function, a piecewise linear function and the maximum over piecewise linear functions. 
    With this result, we can directly apply Theorem~\ref{theo:kl}, using $l(\mathbf{w}) := \lVert\mathbf{w} \rVert^2$ and $\mathbf{A} = \mathbf{0} \in \mathbb{R}^{m \times n}$ such that $l(\mathbf{A} \mathbf{w}) = l(\mathbf{0}) = 0$.  
    This shows that any piecewise linear function is a \gls{KL} function with exponent $\theta = \frac{1}{2}$, and it follows that $\phi$ is a \gls{KL} function with exponent $\theta_{\phi} = \frac{1}{2}$. Further, $\phi + l_\Delta$ is also a \gls{KL} function with exponent $\frac{1}{2}$ [cf. \citet[Section 5]{li2018}]. 

    For the rest of the proof, following the standard techniques developed in \citet{attouch2010}, \citet{aragon2018}, \citet{li2018} and \citet{aragon2020}, the linear convergence of $\{\mathbf{w}_k\}$ to a \gls{KKT} point can be established. 
    We skip the detailed proof for brevity.
\end{prf}


\begin{rema}{rema:kkt_point_local_minimum}
    If $h$ is differentiable at $\bar{\mathbf{w}}$, then the \gls{KKT} point is a local minimizer for\gls{DC} problems with polyhedral components [cf. \citet[p. 4]{PhamDinh2014}]. 
    In our problem, the individual components $g$ and $h$ are no longer piecewise linear due to the regularization term~$\psi$.
    However, the \gls{KKT} system in the limit shown in \eqref{eq:KKT_Cond3}, corresponds to the \gls{KKT} system of the polyhedral \gls{DC} problem without the regularization term, since $\nabla\psi$ vanishes in the limit. 
    Since $h$ is convex, it is differentiable everywhere except on a set of measure zero [cf. \citet[p. 4]{PhamDinh2014}].  
    Thus, we may also conclude that $\bar{\mathbf{w}}$ is almost surely a local minimizer of~\eqref{eq:var_dc}.
\end{rema}


\begin{rema}[Feasibility of Assumptions for Portfolio Optimization Problem]{rema:feasibility_of_assumptions}
    To establish the linear convergence rate, it suffices that either $\nabla g$ \textit{or} $\nabla h$ is locally Lipschitz continuous in a neighborhood of~$\bar{\mathbf{w}}$. 
    For the optimization problem in \eqref{eq:var_dc},  this condition is not always satisfied globally, as both $g$ and $h$ may fail to be continuously differentiable for all $\mathbf{w} \in \Delta$. 
    Nevertheless, this remark demonstrates that under specific structural assumptions, the condition is always met locally. 
    For this purpose, we assume that there are no ties around the $(k +1)^{th}$ scenario, i.e., the portfolio returns are strictly increasing $\bar{\mathbf{w}}^\top \mathbf{x}_k < \bar{\mathbf{w}}^\top \mathbf{x}_{k+1} < \bar{\mathbf{w}}^\top \mathbf{x}_{k+2}$ in this neighborhood.  
    Under this assumption, the discrete $\text{\gls{VaR}}_{X, \alpha}$ is uniquely defined, and we have 
    \begin{equation}\label{eq:cont_diff_CVaR}
        \partial \text{\gls{CVaR}}_{X, \alpha} (\bar{\mathbf{w}}) = \Bigl \{ \nabla \text{\gls{CVaR}}_{X, \alpha} (\bar{\mathbf{w}}) \Bigr\} = \Bigl \{ \frac{1}{\alpha} \cdot \sum_{i=1}^k \mathbf{x}_i \cdot p_i + \frac{\varepsilon}{\alpha} \cdot \mathbf{x}_{k + 1} \Bigr\},
    \end{equation}
     since in this case, the mapping $\bar{\mathbf{w}} \rightarrow \text{\gls{CVaR}}_{X, \alpha} (\bar{\mathbf{w}}) $ is linear in a neighborhood of~$\bar{\mathbf{w}}$ [cf. \citet[p. 871]{wozabal2012}]. 
     Based on this result, we obtain
     \begin{equation}
         \partial h (\bar{\mathbf{w}}) = \Bigl \{  \nabla h(\bar{\mathbf{w}}) \Bigr\} = \Bigl \{ \nabla \psi (\bar{\mathbf{w}}) - \frac{\tau}{\gamma} \cdot \sum_{i=1}^k \mathbf{x}_i \cdot p_i + \frac{\tau}{\gamma} \cdot \mathbf{x}_{k + 1} \cdot (\varepsilon - \gamma)\Bigr\},
     \end{equation}
     which confirms that $h$ is continuously differentiable in a neighborhood of~$\bar{\mathbf{w}}$ [cf. \citet[p. 871]{wozabal2012}]. 
     Since $\nabla h$ is linear in this neighborhood, it follows that it is locally Lipschitz continuous around~$\bar{\mathbf{w}}$. 
     
     \noindent A similar result can be obtained for $\nabla g$, but it requires an additional assumption. 
     In general, the sub-differential of $g$ at $\bar{\mathbf{w}}$ can be expressed as    
     \begin{equation}
         \partial g(\bar{\mathbf{w}}) = \nabla \psi (\bar{\mathbf{w}}) + \nabla f(\bar{\mathbf{w}}) - \tau \cdot \partial \zeta (\bar{\mathbf{w}})
     \end{equation}
    where 
     \begin{equation*}
         \partial \zeta (\mathbf{w}) :=
         \begin{cases}
             \frac{\alpha}{ \gamma} \cdot \nabla \text{\gls{CVaR}}_{X, \alpha} (\bar{\mathbf{w}}) \ &\text{if} \ \frac{\alpha}{ \gamma} \cdot \text{\gls{CVaR}}_{X, \alpha} (\bar{\mathbf{w}}) - r_{\text{min}} < \frac{\alpha - \gamma}{\gamma} \cdot  \text{\gls{CVaR}}_{X, \alpha - \gamma} (\bar{\mathbf{w}}) \\[1em]
             \frac{\alpha - \gamma}{\gamma} \cdot \nabla \text{\gls{CVaR}}_{X, \alpha - \gamma} (\bar{\mathbf{w}}) \ &\text{if} \  \frac{\alpha}{ \gamma} \cdot \text{\gls{CVaR}}_{X, \alpha} (\bar{\mathbf{w}}) - r_{\text{min}} > \frac{\alpha - \gamma}{\gamma} \cdot  \text{\gls{CVaR}}_{X, \alpha - \gamma} (\bar{\mathbf{w}})  \\[1em]
             H_{\bar{\mathbf{w}}} &\text{if} \ \frac{\alpha}{ \gamma} \cdot \text{\gls{CVaR}}_{X, \alpha} (\bar{\mathbf{w}}) - r_{\text{min}} =  \frac{\alpha - \gamma}{\gamma} \cdot  \text{\gls{CVaR}}_{X, \alpha - \gamma} (\bar{\mathbf{w}}),
         \end{cases}
     \end{equation*}
    with $ H_{\bar{\mathbf{w}}} := \text{conv} \Bigl ( \bigl \{ \frac{\alpha}{ \gamma} \cdot \nabla \text{\gls{CVaR}}_{X, \alpha} (\bar{\mathbf{w}}), \frac{\alpha - \gamma}{\gamma} \cdot \nabla \text{\gls{CVaR}}_{X, \alpha - \gamma} (\bar{\mathbf{w}}) \bigr \} \Bigr)$.  
    Hence, $g$ is continuously differentiable for every $\bar{\mathbf{w}}$ where $\frac{\alpha}{ \gamma} \cdot \text{\gls{CVaR}}_{X, \alpha} (\bar{\mathbf{w}}) - r_{\text{min}} \neq \frac{\alpha - \gamma}{\gamma} \cdot  \text{\gls{CVaR}}_{X, \alpha - \gamma} (\bar{\mathbf{w}})$. 
    To make this condition more explicit, consider the case where equality holds 
    \begin{equation}
        \frac{\alpha}{ \gamma} \cdot \text{\gls{CVaR}}_{X, \alpha} (\bar{\mathbf{w}}) - r_{\text{min}} =  \frac{\alpha - \gamma}{\gamma} \cdot  \text{\gls{CVaR}}_{X, \alpha - \gamma} (\bar{\mathbf{w}}),
    \end{equation}
    which, using the discrete representation of the $\text{\gls{CVaR}}_{X, \alpha}$, can be rewritten as
    \begin{equation}
            \biggl ( \frac{1}{\gamma} \cdot \sum_{j=1}^{k} \mathbf{w}^\top \mathbf{x}_{j} \cdot p_{j} + \frac{\varepsilon}{\gamma} \cdot \mathbf{w}^\top \mathbf{x}_{k + 1}  \biggr) - r_{\text{min}} = \biggl( \frac{1}{\gamma} \cdot \sum_{j=1}^{k} \mathbf{w}^\top \mathbf{x}_{j} \cdot p_{j} + \frac{\varepsilon - \gamma}{\gamma} \cdot \mathbf{w}^\top \mathbf{x}_{k + 1} \biggr), 
    \end{equation}
    and further simplified to 
    \begin{equation}
            \frac{\varepsilon}{\gamma} \cdot \mathbf{w}^\top \mathbf{x}_{k + 1}   - r_{\text{min}} =  \frac{\varepsilon - \gamma}{\gamma} \cdot \mathbf{w}^\top \mathbf{x}_{k + 1}.  
    \end{equation}
    This yields the condition 
    \begin{equation}
        \bar{\mathbf{w}}^\top \mathbf{x}_{k + 1} \neq r_{\text{min}},
    \end{equation}
     i.e., $g$ is continuously differentiable at~$\bar{\mathbf{w}}$ whenever the discrete $\text{\gls{VaR}}_{X, \alpha}$ at $\bar{\mathbf{w}}$ is not exactly equal to the lower bound $r_{\text{min}}$. 
     Together with the condition established in \eqref{eq:cont_diff_CVaR}, we obtain that $\nabla g$ is linear in the neighborhood of~$\bar{\mathbf{w}}$, and thus, $\nabla g$ is locally Lipschitz continuous around~$\bar{\mathbf{w}}$.
\end{rema}


\begin{rema}[Generalization of Convergence Result]{rema:generalization-of_convergence_result}
    The results of Theorem~\ref{theo:main1} and Theorem \ref{theo:linear_rate} also apply to problems with piecewise linear objective functions defined over a compact set with linear equality and inequality constraints that fulfill the following requirements: 
    \begin{enumerate}
        \item[(i)] Both functions $g$ and $h$ are strongly convex with the same modulus $\rho > 0$.
        \item[(ii)] The functions $h$ is sub-differentiable at every point in  $\mathrm{dom}$ $h$, i.e. $\partial h(\mathbf{w}) \neq \emptyset$ for all $\mathbf{w} \in \text{dom} \ h$.
        \item[(iii)]    The function $g$ is sub-differentiable at every point in  $\mathrm{dom}$ $g$, i.e. $\partial g(\mathbf{w}) \neq \emptyset$ for all $\mathbf{w} \in \text{dom} \ g$.
        \item[(iv)] The function $g$ or the function $h$ is continuously differentiable at the limit point $\bar{\mathbf{w}}$ with locally Lipschitz derivative.
        \item[(v)] The objective function is bounded from below, i.e., $\underset{\mathbf{w} \in \mathbb{R}^n}{\text{inf}} \phi(\mathbf{w}) > - \infty$.
        \item[(vi)] All feasible points fulfill the \gls{MFCQ}.
    \end{enumerate}
\end{rema}

\newpage


\section{Practical Application \& Numerical Experiments}\label{sec:prac_app}



Having examined the theoretical foundations, the focus now shifts to the practical application. 
By integrating a line search procedure, the \gls{BDCA} is expected to achieve faster convergence and greater robustness in handling complex constraint settings, compared to the \gls{DCA}.
So far, no empirical evaluation has been conducted to assess performance differences between these two \gls{DC} approaches for the optimization problem described in Section~\ref{sec:var_as_dc}.
Likewise, neither method has yet been compared with solvers from chance-constrained programming. 
To address these research gaps, this section presents a structured performance evaluation across multiple data sets, with a focus on challenging constraint settings.


To provide a fair and unbiased algorithmic comparison, the numerical experiments were designed according to the recommendations of \citet{beiranvand2017}. 
In line with these best practices, performance is assessed across three dimensions: efficiency, reliability, and output quality. 
Efficiency is measured via computational time, while output quality is measured by the distance to the best known objective function value. 
Note that in the theoretical part of this paper, solution quality is characterized by convergence to a \gls{KKT} point, which is guaranteed once $\lVert \mathbf{d}_k \rVert$ approaches zero.
In our practical implementation, convergence is monitored using different stopping criteria that assess changes in $\lVert \mathbf{d}_k \rVert$ either directly, or indirectly via in the sequences $\{\mathbf{w}_k\}$ or $\{\phi \left (\mathbf{w}_k \right)\}$, respectively.
Once the algorithm has converged, the focus is on lower objective function values, since for the practical problem at hand, the key performance measure is typically the expected portfolio profit [cf. \citet[p. 155]{Palomar2025}].
Analogously, reliability is evaluated through multiple runs with different starting points, reporting objective function values and the proportion of feasible solutions with 95\% \glspl{CI} [cf. \citet[pp. 655--656]{fox2015}; \citet{Becher1993}; \citet{efron1979}]. 
These intervals are constructed using 100,000 bootstrap samples, and are adjusted using the Bonferroni correction [cf. \citet{Armstrong2014}; \citet{bonferroni1936}].


Since initialization can significantly influence optimization outcomes, the starting points for the reliability evaluation must be chosen carefully. 
Within the upcoming analyses, we experimented with random draws from a Dirichlet distribution to cover different areas of the simplex. 
In general, two initialization schemes are considered: (i) nearly uniform weights across all investment positions, and (ii) highly skewed weights, where only a few decision variables are nonzero. 
The latter setup increases the likelihood of unfavorable initializations, which provides a more rigorous test of the algorithm’s robustness. 
Moreover, we deliberately allow for initializations that may violate the $\text{\gls{VaR}}_{X, \alpha}$ constraint, as in highly restrictive settings, it may not be guaranteed that a feasible starting point can be easily generated.
Table~\ref{tab:exemp_weights} illustrates five examples per initialization scheme for a case with six decision variables.


\begin{table}[ht]
    \scriptsize
    \centering
    \captionsetup{justification=centering}
    \caption{Five random examples for the starting point per initialization scheme.}
    \label{tab:exemp_weights}
    \begin{tabular}{cccccccc}
        \toprule
        Scheme & Draw & $w_{1}$ & $w_{2}$ & $w_{3}$ & $w_{4}$ & $w_{5}$ & $w_{6}$ \\
        \midrule
        \multirow[c]{5}{*}{1} & 1 & 0.17493 & 0.16305 & 0.16348 & 0.16067 & 0.17559 & 0.16228 \\
         & 2 & 0.16894 & 0.17088 & 0.16166 & 0.16667 & 0.16552 & 0.16631 \\
         & 3 & 0.17525 & 0.16799 & 0.16425 & 0.16384 & 0.16547 & 0.16321 \\
         & 4 & 0.16845 & 0.17085 & 0.16294 & 0.17190 & 0.16597 & 0.15989 \\
         & 5 & 0.16823 & 0.16418 & 0.17899 & 0.16459 & 0.16118 & 0.16282 \\
        \midrule
        \multirow[c]{5}{*}{2} & 1 & 0.37573 & 0.00000 & 0.00001 & 0.00012 & 0.22837 & 0.39577 \\
         & 2 & 0.02160 & 0.21907 & 0.01499 & 0.00001 & 0.00051 & 0.74382 \\
         & 3 & 0.00791 & 0.00001 & 0.99208 & 0.00000 & 0.00000 & 0.00000 \\
         & 4 & 0.39894 & 0.47972 & 0.00978 & 0.00000 & 0.02959 & 0.08197 \\
         & 5 & 0.00000 & 0.00000 & 0.01100 & 0.98892 & 0.00008 & 0.00000 \\
        \bottomrule
    \end{tabular}
\end{table}


To validate the general applicability of the proposed methods, the numerical experiments incorporate distinct real-world data sets. 
These data sets were originally compiled and preprocessed by \citet{BRUNI2016}, and include adjustments for dividends and stock splits.
Furthermore, their public availability facilitates reproducibility and enables future benchmarking of results.
Table~\ref{tab:data_descrip} provides a brief overview of the data. 
Notably, two of these data sets contain fewer decision variables but a larger number of scenarios, whereas the remaining two comprise a broader set of index positions, with comparatively fewer scenarios.
The data set names refer to collections of individual stocks in the case of \gls{DJ}, \gls{NASDAQ}, and \gls{FTSE}, and to industry-based portfolios in the case of \gls{FF49}.
More information about the data derivation can be found in \citet{BRUNI2016}.


\begin{table}[H]
    \scriptsize
    \centering
    \captionsetup{justification=centering}
    \caption{Description of Real-World Data Sets.  \\ \emph{Source:} cf. \citet[p. 860]{BRUNI2016}.}
    \label{tab:data_descrip}
        \begin{tabular}{lccrcr}
            \toprule
            Data         & Country   & n     & \multicolumn{1}{c}{S}        & Interval  & \multicolumn{1}{c}{Time Period}  \\  
            \midrule
            
            \gls{DJ}        & USA       & 28    & 1,363    & Weekly    & Feb 1990  - Apr 2016     \\ 
            \gls{FF49}      & USA       & 49    & 2,325    & Weekly    & Jul 1969 - Apr 2016      \\ 
            \gls{FTSE}      & UK        & 83    & 717      & Weekly    & Jul 2002 - Apr 2016    \\               \gls{NASDAQ}    & USA       & 82    & 596      & Weekly    & Nov 2004 - Apr 2016    \\ 

            \bottomrule
        \end{tabular}
\end{table}


To obtain more detailed insights into the data, Table \ref{tab:stats} provides summary statistics for each index. 
Note that the mean $(\mu)$, standard deviation $(\sigma)$, $\text{\gls{VaR}}_{0.05}$,  $\text{\gls{CVaR}}_{0.05}$, skewness $(\gamma_1)$ and kurtosis $(\gamma_2)$ have been calculated individually per index position. 
Table \ref{tab:stats} summarizes the distribution of these metrics across all positions by reporting the minimum, maximum, as well as the 25th, 50th, and 75th percentiles (Q25, Q50, Q75). 
Overall, the results for kurtosis indicate that the distributions of the investment positions, regardless of the index, exhibit fatter tails compared to a normal distribution. 
This implies a higher likelihood of extreme gains or losses.
Similarly, the summary statistics for skewness show that many distributions are at least moderately asymmetric. 
This suggests that certain index positions provide higher chances for either large gains or large losses.
Particularly pronounced positive and negative skewness can be observed in the \gls{FTSE} and \gls{NASDAQ} data. 
These two indices also exhibit the largest risk variety, as the $\text{\gls{VaR}}_{0.05}$ and $\text{\gls{CVaR}}_{0.05}$ distributions highlight several positions with elevated downside risk.
All in all, the diversity in statistical profiles across the data sets creates a more complex evaluation environment, such that different algorithmic properties may emerge throughout the comparisons. 


\begin{table}[H]
    \scriptsize
    \centering
    \captionsetup{justification=centering}
    \caption{Summary Statistics of Real-World Data Sets.}
    \label{tab:stats}

    \begin{tabular}{cccccccc}
    \toprule
     Data & Stats & $\mu$ & $\sigma$ & $\text{\gls{VaR}}_{0.05}$ & $\text{\gls{CVaR}}_{0.05}$ & $\gamma_1$ & $\gamma_2$ \\
    \midrule
    \multirow[c]{5}{*}{\gls{DJ}} & Min & 1.00130 & 0.02828 & 0.91010 & 0.86600 & -1.37390 & 4.39180 \\
     & Q25 & 1.00210 & 0.03317 & 0.93310 & 0.90150 & -0.09990 & 5.64510 \\
     & Q50 & 1.00250 & 0.04000 & 0.94340 & 0.91440 & 0.08350 & 6.50290 \\
     & Q75 & 1.00280 & 0.04583 & 0.95090 & 0.92900 & 0.24140 & 7.37100 \\
     & Max & 1.00610 & 0.06164 & 0.95880 & 0.93720 & 0.74870 & 18.07420 \\
    \midrule
    \multirow[c]{5}{*}{\gls{FF49}} & Min & 1.00280 & 0.01732 & 0.92470 & 0.88230 & -0.63390 & 6.01520 \\
     & Q25 & 1.00370 & 0.02828 & 0.95460 & 0.92750 & -0.20400 & 8.47560 \\
     & Q50 & 1.00390 & 0.03000 & 0.95780 & 0.93450 & -0.07680 & 11.07510 \\
     & Q75 & 1.00440 & 0.03317 & 0.96150 & 0.93900 & 0.25470 & 15.09310 \\
     & Max & 1.00540 & 0.05477 & 0.97630 & 0.95930 & 1.96890 & 37.92980 \\
    \midrule
    \multirow[c]{5}{*}{\gls{FTSE}} & Min & 0.99940 & 0.02449 & 0.89550 & 0.82070 & -1.86500 & 4.18180 \\
     & Q25 & 1.00150 & 0.03606 & 0.92210 & 0.88160 & -0.10900 & 5.77730 \\
     & Q50 & 1.00250 & 0.04000 & 0.93830 & 0.91010 & 0.13930 & 7.19970 \\
     & Q75 & 1.00360 & 0.05477 & 0.94780 & 0.92430 & 0.37110 & 10.25650 \\
     & Max & 1.00800 & 0.08775 & 0.96210 & 0.94570 & 2.23860 & 27.53910 \\
    \midrule
    \multirow[c]{5}{*}{\gls{NASDAQ}} & Min & 1.00010 & 0.02646 & 0.88360 & 0.82100 & -1.16480 & 4.07070 \\
     & Q25 & 1.00190 & 0.03873 & 0.92200 & 0.88170 & -0.02490 & 5.22690 \\
     & Q50 & 1.00310 & 0.04690 & 0.93170 & 0.89850 & 0.17690 & 6.41960 \\
     & Q75 & 1.00470 & 0.05385 & 0.94200 & 0.91530 & 0.49860 & 8.27620 \\
     & Max & 1.01030 & 0.08367 & 0.96360 & 0.94450 & 1.55850 & 22.79620 \\
    \bottomrule
    \end{tabular}

\end{table}


Building on the theoretical framework outlined above, the following three subsections present numerical experiments that individually examine: (i) the influence of input parameters on the performance of \gls{BDCA}, (ii) comparative results with respect to \gls{DCA}, and (iii) benchmark comparisons with solvers from chance-constrained programming.
The computational experiments were carried out using the \href{https://www.southampton.ac.uk/research/facilities/iridis-research-computing-facility}{IRIDIS~5 Research Computing Facility}, where each compute node is equipped with an Intel Xeon Gold 6130 processor (64 cores at 2.1 GHz) and DDR4 memory. 
The algorithms were implemented in Python~3.12.7 using the \texttt{CVXPY} optimization framework provided by \citet{diamond2016}. 
To facilitate reproducibility, the complete implementation is publicly available via the project’s GitHub repository:
\href{https://github.com/mlthormann/BDCA-For-Portfolio-Optimization/}{BDCA-For-Portfolio-Optimization \faExternalLink}.



\subsection{Sensitivity Analysis of BDCA Input Parameters}

\begin{figure}[H]
    \centering
    \captionsetup{justification=centering}
    \includegraphics[width=1\textwidth]{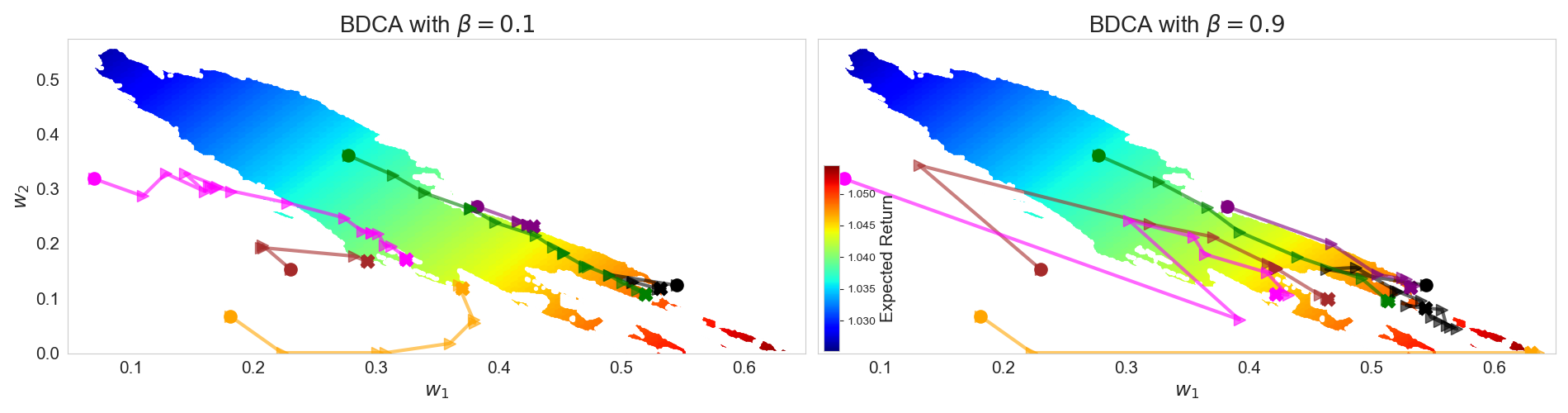}
    \caption{Synthetic data example to showcase the influence of $\beta$, with $\rho = 0.5$ and $\tau = 1.5$ held constant.}%
    \label{fig:influence_of_beta}
\end{figure}

Before comparing \gls{BDCA} with other methods, this subsection investigates how its performance depends on the choice of input parameters $\tau$, $\rho$, and $\beta$.
The overall goal is to identify configurations that yield robust results -- even in the presence of challenging $\text{\gls{VaR}}_{X, 0.05}$ constraints or under adverse initialization.
For this purpose, we begin with a concise theoretical overview of how each parameter affects the algorithmic robustness.
Among the three inputs, $\beta$ controls the reduction of the extrapolation step size~$\lambda_k$ during the line search procedure [cf. Algorithm~\ref{alg:bdca}].
Larger values promote more aggressive extrapolation, which can increase the number of function evaluations but may also accelerate convergence. Moreover, it allows for a broader exploration  of the solution space, such that feasible outcomes may be found even in challenging settings.
In contrast, smaller values reduce the number of function evaluations but can slow down progress due to more conservative step sizes.
Figure~\ref{fig:influence_of_beta} illustrates the impact of $\beta \in \{0.1, 0.9\}$ on a synthetic three-variable portfolio problem, where the shaded region denotes the feasible set under the selected $\text{\gls{VaR}}_{X, 0.05}$ constraint.
In this toy example, the trajectories for multiple initializations demonstrate that a large value for $\beta$ leads to faster convergence and broader exploration of the feasible region. 


\vspace{-0.2cm}

\begin{figure}[H]
    \centering
    \captionsetup{justification=centering}
    \includegraphics[width=1\textwidth]{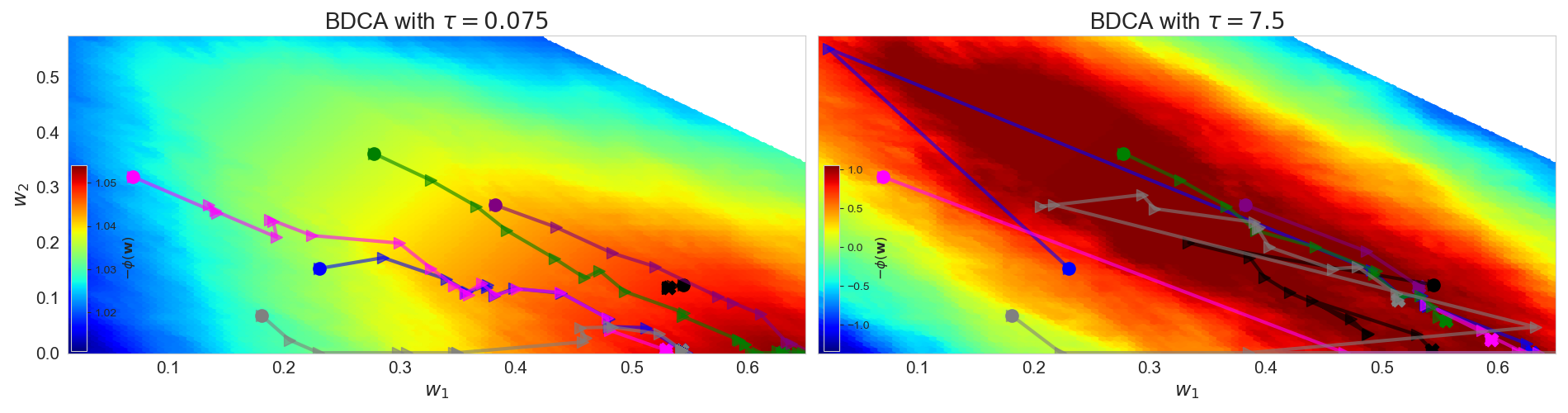}
    \caption{Synthetic data example to showcase the influence of $\tau$, with $\beta = 0.9$ and $\rho = 0.5$ held constant.}%
    \label{fig:influence_of_tau}
\end{figure}

Next, we consider the penalty parameter, which incorporates the $\text{\gls{VaR}}_{X, \alpha}$ constraint into the objective function shown in \eqref{eq:penalized_problem}.
Unlike $\beta$, $\tau$ impacts both the convex subproblem as well as the line search procedure. 
Larger values emphasize constraint satisfaction but can also disproportionately influence the objective function. 
Specifically, the penalty term can vanish abruptly once feasibility is achieved, or minor improvements in the constraint satisfaction can dominate the line search, potentially  causing erratic descent behavior.
Smaller values mitigate this effect by creating a smoother penalization landscape and reducing sensitivity to constraint violations.
However, they can also slow down convergence and feasibility may not be guaranteed under challenging settings.
Figure~\ref{fig:influence_of_tau} illustrates the effect of $\tau \in \{0.075, 7.5\}$ using the same synthetic portfolio setup without excluding infeasible solutions.
The right-hand panel demonstrates how a larger penalty parameter emphasizes the feasible region but can also lead to unstable algorithmic behavior.
Conversely, the left-hand panel shows how a too small $\tau$ prolongs progress toward feasibility.


\vspace{-0.2cm}

\begin{figure}[H]
    \centering
    \captionsetup{justification=centering}
    \includegraphics[width=1\textwidth]{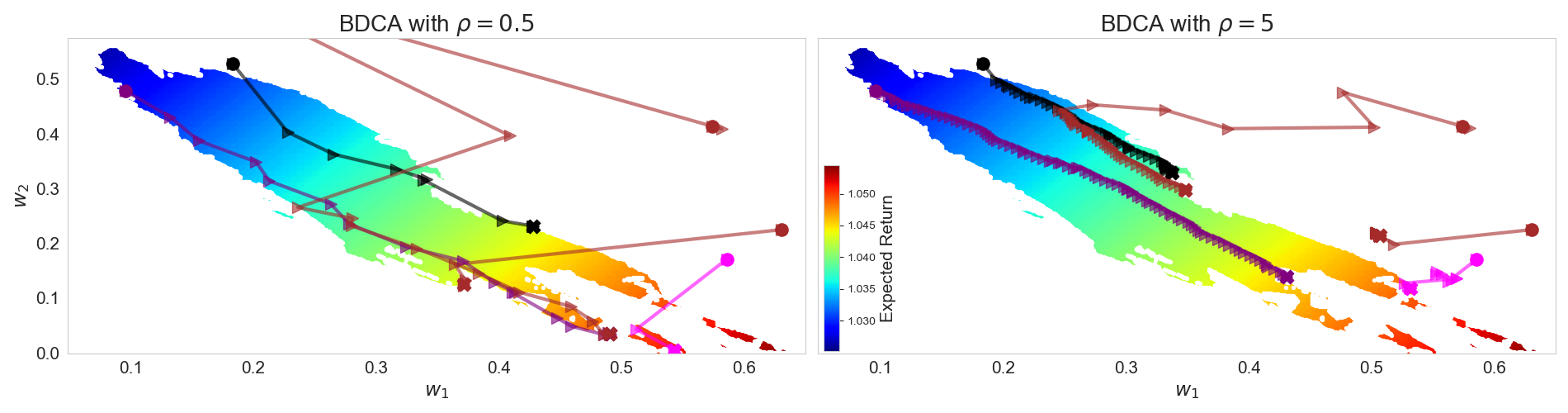}
    \caption{Synthetic data example to showcase the influence of $\rho$, with $\beta = 0.75$ and $\tau = 1.5$ held constant.}%
    \label{fig:influence_of_rho}
\end{figure}

Similar to $\tau$, the parameter $\rho$ plays a dual role in the optimization process.
It is included in the regularization term $\psi(\mathbf{w})$, which appears in both \gls{DC} components $g(\mathbf{w})$ and $h(\mathbf{w})$. 
Although $\psi (\mathbf{w})$ cancels out in the overall objective, it still influences the convex subproblem~\eqref{eq:convex_subproblem}, where it penalizes extreme portfolio weights and encourages proximity between successive iterates.
Moreover, $\rho$ affects the Armijo-type condition in the line search, where high values make it more difficult to terminate the backtracking procedure.
Overall, larger values of $\rho$ promote smoother updates and enhance numerical stability, but excessive regularization can suppress meaningful variation and slow down convergence.
This trade-off is captured in Figure~\ref{fig:influence_of_rho} using the same synthetic data example with $\rho \in \{0.5, 5\}$.
The right-hand panel shows that strong regularization yields stable trajectories but limits the progress toward areas with higher portfolio returns and feasibility is not guaranteed for all outcomes.
In contrast, the left-hand panel demonstrates that a smaller $\rho$ allows more variation in the optimization path but can lead to erratic steps.


\vspace{-0.2cm}

\begin{figure}[H]
    \centering
    \includegraphics[width=1\textwidth]{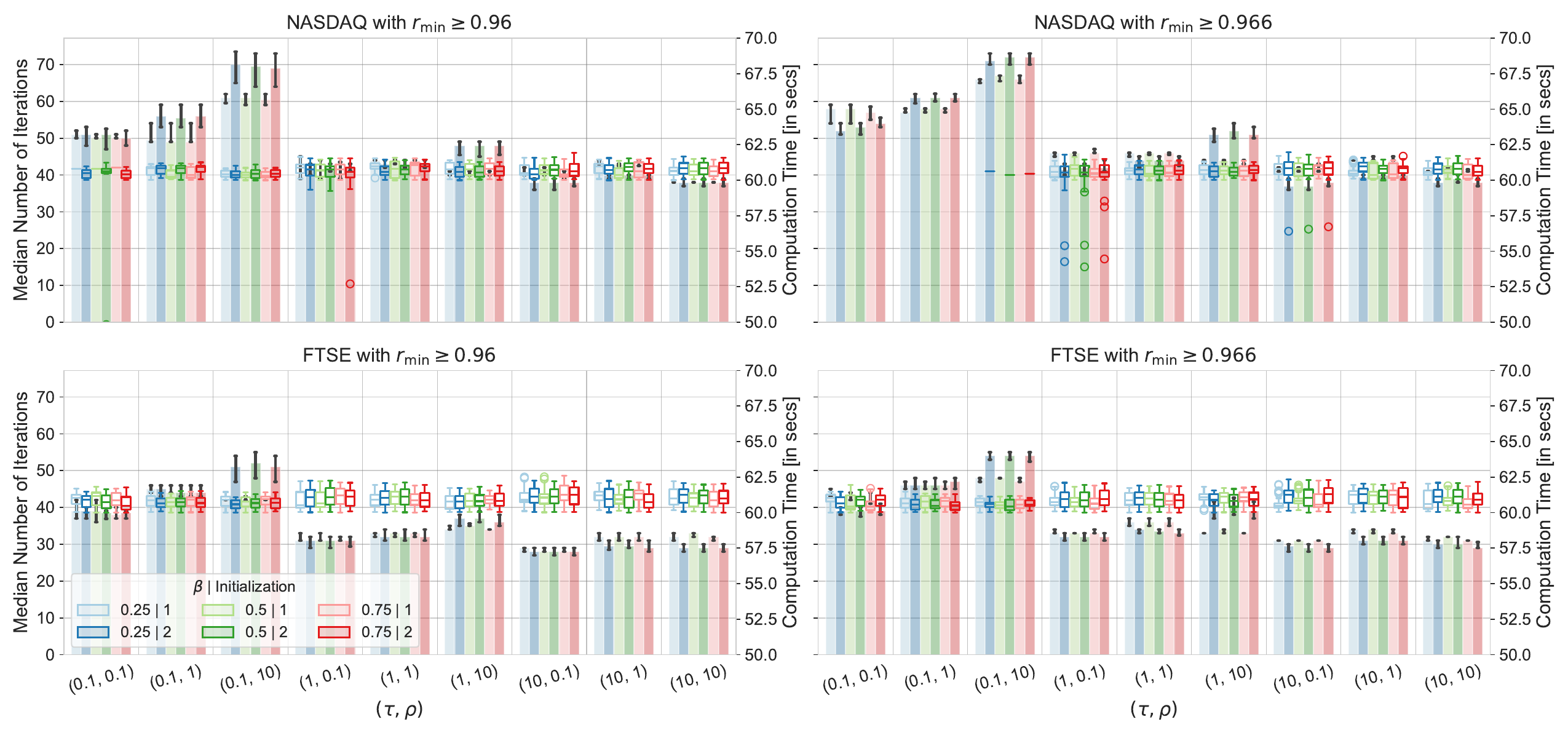}
    \caption{Impact of \gls{BDCA} parameters on the computation time and iteration count for the \gls{NASDAQ} and \gls{FTSE} data. 
    Each subplot displays the median number of iterations (left y-axis, bar plots) and the computation time of all feasible outcomes (right y-axis, box plots) across varying combinations of $\tau$ and $\rho$ on the x-axis. 
    Results are grouped by $\beta \in \{0.25, 0.5, 0.75\}$ and initialization strategy (1 or 2). }%
    \label{fig:sensitivity_analysis_bdca_2}
\end{figure}

Building on the theoretical considerations, we perform a sensitivity analysis to evaluate the performance of \gls{BDCA} across all combinations of the parameters $\beta \in \{0.25, 0.5, 0.75\}$, $\rho \in \{0.1, 1, 10\}$, and $\tau \in \{0.1, 1, 10\}$, using the four data sets introduced earlier.
Each configuration is tested with 50 randomly generated starting points per initialization scheme. 
To manage computational cost, runs are terminated after 60 seconds if convergence is not reached.
Figure~\ref{fig:sensitivity_analysis_bdca_2} presents computation times of all feasible solutions (right y-axis, box plots) and iteration counts (left y-axis, bar plots) for the \gls{NASDAQ} and \gls{FTSE} data across all configurations.
In most cases, the algorithm reached the 60-second time limit, which results in consistently high computation times with very similar variability. 
The iteration counts reflect that smaller values of $\tau$ result in more iterations within the same time frame, particularly when combined with a higher $\rho$.
This may be due to numerical difficulties encountered by the convex solver of the subproblem when $\tau$ is large.
Corresponding results for the remaining data sets are shown in Figure~\ref{fig:sensitivity_analysis_bdca_2_dj_ff49} in Appendix~\ref{app:numerical_results}.


\vspace{-0.2cm}

\begin{figure}[H]
    \centering
    \includegraphics[width=1\textwidth]{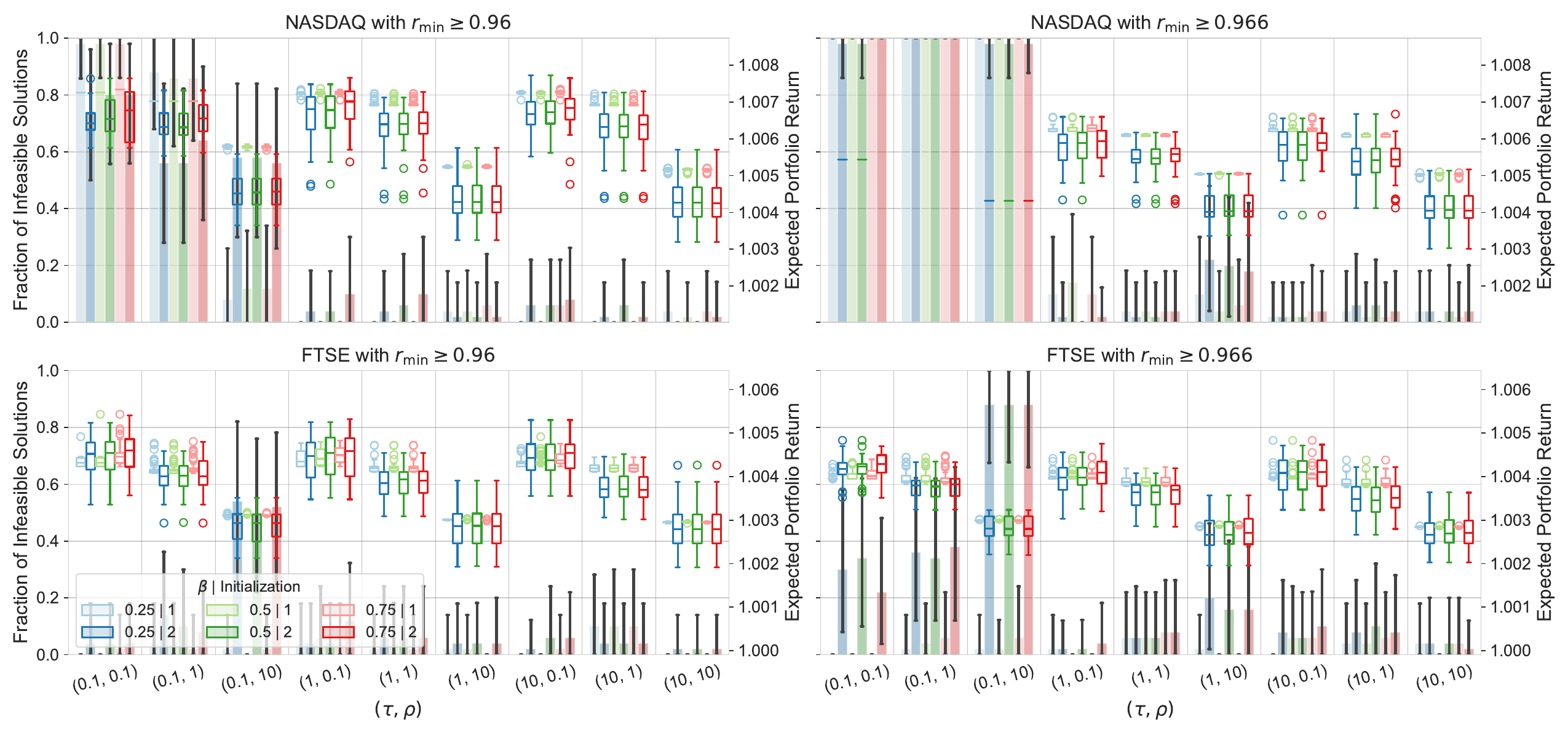}
    \caption{Impact of \gls{BDCA} parameters on feasibility and expected return for the \gls{NASDAQ} and \gls{FTSE} data. 
    Each subplot displays the fraction of infeasible solutions (left y-axis, bar plots) and the expected portfolio returns of all feasible outcomes (right y-axis, box plots) across varying combinations of $\tau$ and $\rho$ on the x-axis. 
    Results are grouped by $\beta \in \{0.25, 0.5, 0.75\}$ and initialization strategy (1 or 2). }%
    \label{fig:sensitivity_analysis_bdca}
\end{figure}

To complete the sensitivity analysis, Figure~\ref{fig:sensitivity_analysis_bdca} illustrates the fraction of infeasible solutions  (left y-axis, bar plots) and the expected return of all feasible outcomes  (right y-axis, box plots) for the \gls{NASDAQ} and \gls{FTSE} data.
Overall, two main trends can be observed.
First, feasibility is strongly affected by~$\tau$, with $\tau = 0.1$ consistently producing the highest fraction of infeasible solutions.
This effect is predominantly amplified for combinations with larger $\rho$ values and under the second initialization scheme.
Since the penalty parameters is too small the $\text{\gls{VaR}}_{X, 0.05}$ constraint is not sufficiently enforced.
Second, the expected portfolio return is drastically affected by the type of starting point and is more sensitive to $\rho$.
Smaller values generally produce higher returns, as they allow the algorithm to explore more diverse portfolio weights.
In contrast, $\beta$ shows minimal influence, as performance curves across its values are nearly identical.
Analogous results for the remaining data sets are shown in Figure~\ref{fig:sensitivity_analysis_bdca_dj_ff49} in Appendix~\ref{app:numerical_results}.


\vspace{-0.2cm}


The previous results highlight the limitations of fixed parameter choices.
While certain configurations yield acceptable performance on specific instances, no single setting robustly ensures feasibility, fast convergence, and high returns across all data sets and initialization schemes.
To address this limitation, we challenge the performance of the fixed parameter choices with an adaptive mechanism that can respond to problem-specific characteristics during the optimization process. 
Specifically, the penalty parameter is initialized with $\tau_0 = 0.01$, as in most test cases this value does not sufficiently enforce the constraint in the objective. 
Afterwards, the parameter is increased in every iteration where constraint violations persist and the improvement regarding the constraint satisfaction stagnates.
Such procedures are well-established in penalty methods, and ensure that the algorithm identifies a penalty level that is sufficiently large to enforce feasibility but is not larger than necessary  [cf. \citet[ch. 8]{Vassiliadis2021}].
The regularization parameter generally starts at $\rho_0 = 0.01 \cdot n$ and is slightly adjusted based on the extremeness of the starting point to give more space to recover from an adverse initialization. 
Over the iterations, the parameter is then gradually reduced to allow greater flexibility in portfolio weights as optimization progresses.
As the influence of the backtracking parameter was limited in the previous experiments, it is set to $\beta = 0.9$ by default but is adjusted based on logical tests that access the descent quality in each iteration. 
For the exact implementation of the adaptive strategy, the interested reader is referred to the project’s GitHub repository:
\href{https://github.com/mlthormann/BDCA-For-Portfolio-Optimization/}{BDCA-For-Portfolio-Optimization \faExternalLink}.


Table~\ref{tab:perform_ada_strat} summaries the results of the adaptive strategy for all data sets.
Overall, the dynamic parameter adjustments have a  positive impact on all considered performance measures.
A bigger improvement can be observed for the computation times, as most median outcomes are now under the 60-second time limit, independent of the initialization scheme.
Similarly, the number of iterations till convergence could be reduced, especially for the \gls{NASDAQ} examples.
Moreover, with the adaptive strategy, the algorithmic output quality is less dependent on the starting point.
The median returns are higher compared to the results of the static parameters, and the results between initialization schemes are now very similar. 
Lastly, the number of infeasible solutions is very small as most of the outcomes no longer violate the $\text{\gls{VaR}}_{X, 0.05}$ requirement. 
Therefore, the adaptive strategy will be maintained in the upcoming numerical experiments.


\begin{table}[H]
    \scriptsize
    \centering
    \caption{Impact of the adaptive parameter strategy on the performance of \gls{BDCA}. The results are grouped by data set (\gls{FTSE}, \gls{NASDAQ}, \gls{DJ}, \gls{FF49}), $\text{\gls{VaR}}_{X, 0.05}$ threshold (0.96, 0.964) and the initialization strategy (1, 2). Median results with \glspl{CI} are reported for the expected portfolio returns, the iteration counts and the computation times, while the last column shows the number of infeasible (inf.) solutions.}
    \label{tab:perform_ada_strat}
   \begin{tabular}{ccc | ccc | ccc | ccc | c}
    \toprule
     &  &  & \multicolumn{3}{c}{Expected Return} & \multicolumn{3}{c}{Iteration Count} & \multicolumn{3}{c}{Computation Time} & Inf. \\
    Data & $r_{\text{min}}$ & Scheme & Low & Q50 & Up & Low & Q50 & Up & Low & Q50 & Up & Cnt \\
    \midrule
    \multirow[c]{4}{*}{\gls{DJ}} & \multirow[c]{2}{*}{0.96} & 1 & 1.00446 & 1.00451 & 1.00455 & 24 & 27 & 31 & 8.8 & 10.1 & 12.3 & 1 \\
     &  & 2 & 1.00411 & 1.00425 & 1.00442 & 22 & 27 & 32 & 6.9 & 8.7 & 10.7 & 0 \\
    \cmidrule{2-13}
     & \multirow[c]{2}{*}{0.966} & 1 & 1.00368 & 1.00377 & 1.00382 & 28 & 36 & 43 & 7.9 & 10.5 & 13.0 & 1 \\
     &  & 2 & 1.00343 & 1.00356 & 1.0037 & 28 & 41 & 52 & 7.8 & 11.0 & 13.1 & 6 \\
    \midrule
    \multirow[c]{4}{*}{\gls{FF49}} & \multirow[c]{2}{*}{0.96} & 1 & 1.00504 & 1.00529 & 1.00539 & 18 & 28 & 36 & 40.8 & 74.9 & 89.2 & 1 \\
     &  & 2 & 1.00507 & 1.00521 & 1.00533 & 23 & 30 & 47 & 65.6 & 85.0 & 216.5 & 2 \\
    \cmidrule{2-13}
     & \multirow[c]{2}{*}{0.966} & 1 & 1.00484 & 1.00489 & 1.00496 & 24 & 29 & 36 & 69.0 & 88.2 & 108.2 & 0 \\
     &  & 2 & 1.00477 & 1.00487 & 1.00497 & 25 & 31 & 38 & 72.6 & 92.1 & 111.9 & 0 \\
    \midrule
    \multirow[c]{4}{*}{\gls{FTSE}} & \multirow[c]{2}{*}{0.96} & 1 & 1.00535 & 1.00539 & 1.00555 & 18 & 20 & 23 & 25.8 & 28.4 & 34.0 & 2 \\
     &  & 2 & 1.0052 & 1.0054 & 1.00558 & 19 & 21 & 25 & 27.6 & 30.9 & 36.0 & 1 \\
    \cmidrule{2-13}
     & \multirow[c]{2}{*}{0.966} & 1 & 1.00472 & 1.00485 & 1.00504 & 21 & 23 & 27 & 28.5 & 33.3 & 40.5 & 3 \\
     &  & 2 & 1.00451 & 1.00481 & 1.00509 & 23 & 27 & 31 & 32.8 & 39.9 & 49.5 & 5 \\
    \midrule
    \multirow[c]{4}{*}{\gls{NASDAQ}} & \multirow[c]{2}{*}{0.96} & 1 & 1.00751 & 1.00758 & 1.00767 & 22 & 25 & 28 & 26.2 & 30.4 & 33.3 & 0 \\
     &  & 2 & 1.00738 & 1.00765 & 1.00771 & 22 & 25 & 29 & 26.5 & 29.3 & 35.1 & 0 \\
    \cmidrule{2-13}
     & \multirow[c]{2}{*}{0.966} & 1 & 1.00653 & 1.00674 & 1.00684 & 19 & 20 & 24 & 20.1 & 22.2 & 28.3 & 1 \\
     &  & 2 & 1.00619 & 1.00648 & 1.00667 & 26 & 31 & 35 & 29.4 & 37.2 & 43.8 & 3 \\
    \bottomrule
    \end{tabular}
        
\end{table}



\subsection{Comparative Analysis of BDCA and DCA}


After evaluating the sensitivity of \gls{BDCA}'s performance to its input parameters, this subsection examines the overall impact of the line search. 
As noted in the introduction, the \gls{DCA} may fail to identify feasible solutions under challenging constraints, and can become computationally inefficient for large-scale problems. 
To demonstrate that the proposed line search mitigates these limitations, we compare the performance of \gls{BDCA} and \gls{DCA} across the previously introduced data sets under various $\text{\gls{VaR}}_{X, 0.05}$ constraints with increasing levels of difficulty.  
Specifically, we selected the majority of risk thresholds $r_{\text{min}} \in \{0.96, 0.966, 0.968, 0.97, 0.972 \}$ such that for most data sets, the desired $\text{\gls{VaR}}_{X, 0.05}$ cannot be satisfied through an individual index position [cf. Table~\ref{tab:stats}]. 
Throughout the experiments each algorithm is tested with 250 randomly generated starting points per initialization scheme. 
To control computational cost, runs are terminated after 240 seconds if convergence is not achieved. 


Before we proceed to the numerical results, we first discuss how convergence is defined in the rest of this paper. 
\citet{beiranvand2017} emphasize that termination criteria can substantially affect performance comparisons. 
For the first application of the \gls{DCA} to the portfolio selection problem at hand, \citet{wozabal2012} derived finite-step convergence. 
This implies that the normed difference between $\mathbf{w}_k$ and $\mathbf{w}_{k+1}$ will be exactly zero after an unknown number of iterations, i.e., $\lVert \mathbf{d}_k \rVert = 0$.
In contrast, this paper provides a linear convergence rate for applying the \gls{BDCA} to the same optimization problem with an adjusted \gls{DC} decomposition. 
Thus, the generated solution will approach the limit point across iterations, but convergence may not be reached exactly, and an appropriate stopping condition must be chosen. 
Based on our theoretical results, the potential termination criteria can be grouped into two categories. On the one hand, the convergence of the sequence $\{ \mathbf{w}_k \}$ to a stationary point can be verified via 
    \begin{align*}
        \begin{split}
            \text{\gls{dkabs}:}& \quad  \lVert \mathbf{w}_{k + 1} - \mathbf{w}_k \rVert  \leq \epsilon \quad \text{or}\\
            \text{\gls{dkrel}:}& \quad \max_i \bigl \{ \frac{\lvert w_{i, k + 1} - w_{i, k} \rvert}{\lvert w_{i,k} \rvert} \bigr \} \leq \epsilon,
        \end{split}
    \end{align*}
where $\epsilon \in \mathbb{R}_+$ denotes the tolerance that the change between consecutive solution vectors must not exceed. On the other hand, the convergence of the sequence $\{ \phi(\mathbf{w}_k) \}$ to a limit point~$\Bar{\phi}$ could be monitored via
    \begin{align*}
        \begin{split}
            \text{\gls{fctabs}:}& \quad \lVert \phi \left(\mathbf{w}_{k + 1} \right) - \phi \left(\mathbf{w}_{k} \right) \rVert \leq \epsilon \quad \text{or}\\
            \text{\gls{fctcrel}:}& \quad \frac{\lvert \phi\left(\mathbf{w}_{k + 1} \right) - \phi\left(\mathbf{w}_{k}\right) \rvert}{\lvert \phi \left(\mathbf{w}_{k} \right) \rvert} \leq \epsilon, \\
        \end{split}
    \end{align*}
where the change between consecutive objective values must not exceed the tolerance accordingly. 
Within our preliminary analysis, we investigated the impact of these stopping conditions on the three performance criteria, with tolerance set to $\epsilon = 10^{-5}$. 
However, the influence of the stopping conditions on the performance was limited for the selected data sets. 
Therefore, all results in the remaining parts of this paper only show the outputs using \gls{dkabs}, while a detailed comparison between the different stopping criteria is provided in Appendix~\ref{app:numerical_results}.


Figure~\ref{fig:comp_dca_1} compares the fraction of infeasible solutions (left y-axis, bar plots) and the expected portfolio returns of all feasible outcomes (right y-axis, box plots) across the selected constraint settings. 
Overall, the numerical results confirm that the added line search procedure enhances both reliability and output quality. 
The bar plots per initialization scheme and algorithm highlight that the outputs of \gls{DCA} are more sensitive to the starting point: under the more skewed initialization setup, the algorithm systematically finds fewer feasible solutions and yields lower median expected portfolio profits with larger variability.
Compared to the \gls{BDCA}, the \gls{DCA} also struggles significantly to find feasible solutions for the most challenging constraint settings, regardless of the initialization scheme or data set. 
These performance differences are particularly evident for the \gls{DJ} data, but also observable across the other data sets for higher $r_{\text{min}}$ values.
In contrast, the results of \gls{BDCA} indicate that the algorithm finds a similar number of feasible solutions, independent of the starting point. 
Moreover, the median portfolio profits are closer to the best-known solution, and their variability is more consistent across initialization schemes, as the \gls{IQR} between box plots often overlaps.
Compared to the \gls{DCA}, this improved performance can be attributed to larger step sizes enabled by the line search procedure.
These facilitate escaping from infeasible initializations where the $\text{\gls{VaR}}_{X, 0.05}$ constraint is not satisfied and enable the exploration of regions in the solution space with higher returns.


\begin{figure}[H]
    \centering
    \includegraphics[width=1\textwidth]{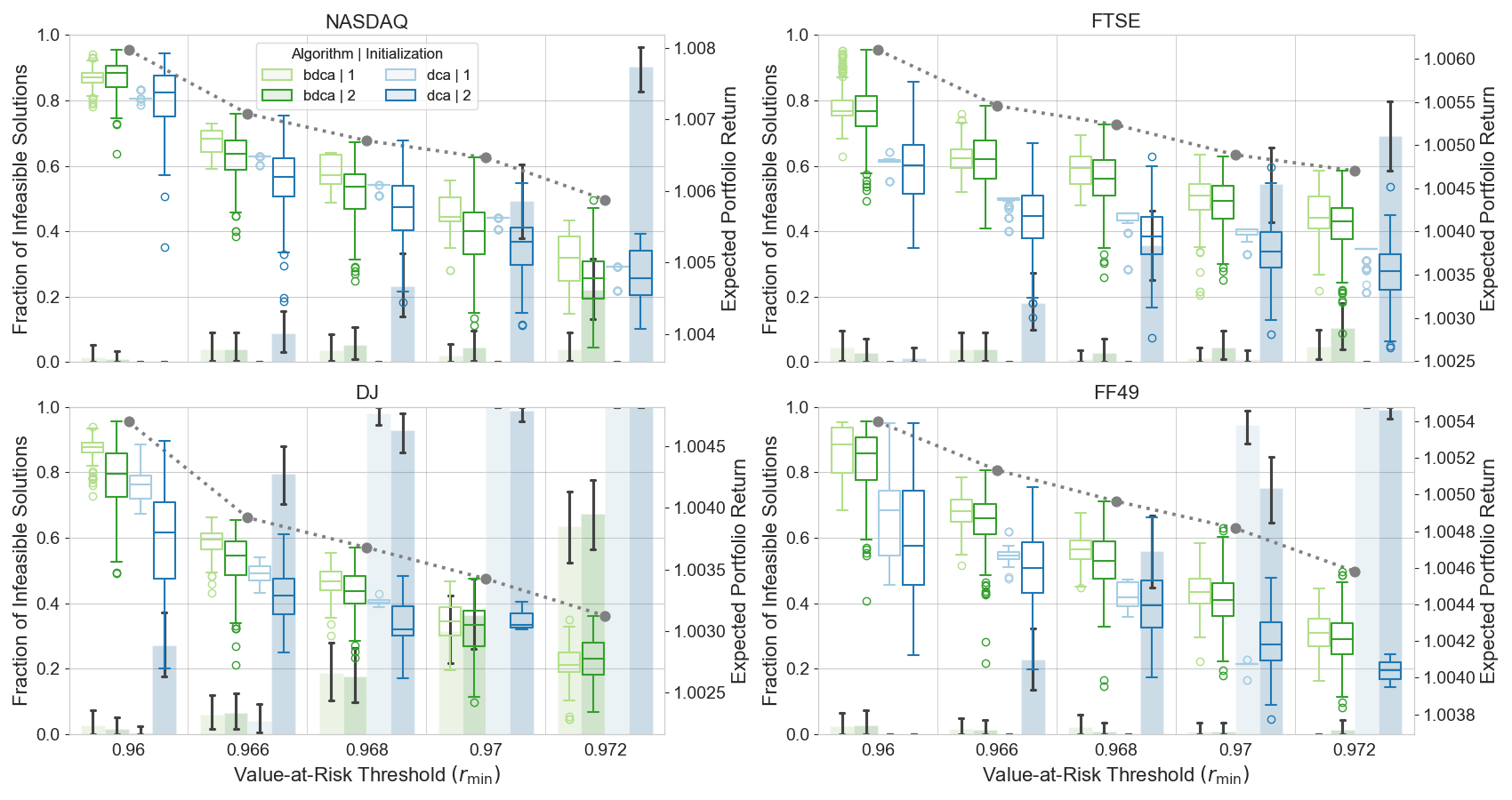}
    \caption{Performance comparison between \gls{BDCA} and \gls{DCA} using financial data sets and focusing on algorithmic reliability and output quality. Each subplot displays the fraction of infeasible solutions (left y-axis, bar plots) and the expected portfolio returns of all feasible outcomes (right y-axis, box plots) across varying $\text{\gls{VaR}}_{X, 0.05}$ levels on the x-axis. Results are grouped by solver type (\gls{BDCA}, \gls{DCA}) and initialization strategy (1 or 2). The gray horizontal line marks the highest portfolio profit achieved among all feasible solutions, independent of the solver.}%
    \label{fig:comp_dca_1}
\end{figure}

\vspace{-0.2cm}

\begin{figure}[H]%
    \centering
    \includegraphics[width=1\textwidth]{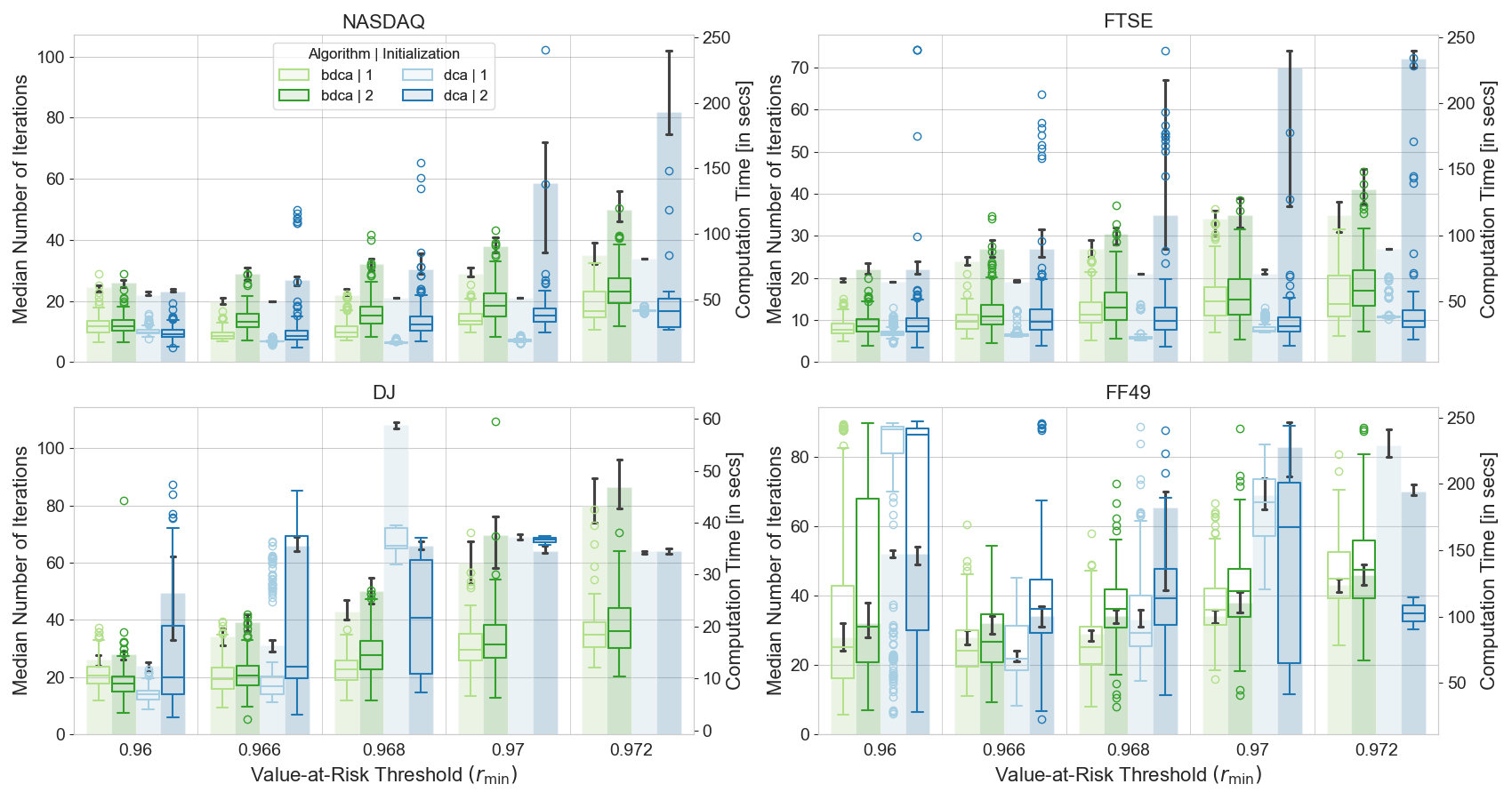}
    \caption{Performance comparison between \gls{BDCA} and \gls{DCA} using financial datasets and focusing on algorithmic efficiency. Each subplot displays the median number of iterations (left y-axis, bar plots) and the computation time of all feasible solutions (right y-axis, box plots) across varying $\text{\gls{VaR}}_{X, 0.05}$ levels (0.96, 0.966, 0.968, 0.97, 0.972) on the x-axis. Results are grouped by solver type (\gls{BDCA}, \gls{DCA}) and initialization strategy (1 or 2).}%
    \label{fig:comp_dca_2}%
\end{figure}

To complete the comparative analysis, Figure~\ref{fig:comp_dca_2} compares the median iteration count (left y-axis, bar plots) and the computation time of all outcomes (right y-axis, box plots) across the selected test cases.
Overall, the numerical results confirm that the added line search procedure also enhances the algorithmic efficiency. 
For the easiest constraint setting and the two data sets with more decision variables, the algorithms exhibit similar computational effort, as run times and iteration counts are comparable.
However, for more demanding $\text{\gls{VaR}}_{X, 0.05}$ levels, the performance variability of the \gls{DCA} increases.
In particular, when combined with the second initialization scheme, the algorithm shows higher median iteration counts and run times, along with a greater number of outliers.
In contrast, the results of the \gls{BDCA} indicate that the line search significantly mitigates performance differences between the initialization schemes.
For the two data sets with fewer investment options but a higher number of scenarios, the described performance discrepancy between the algorithms can be observed across all constraint settings.  
Similar to the results of the algorithmic reliability and output quality, the superiority of the \gls{BDCA} over the \gls{DCA} can be attributed to larger step sizes enabled by the line search procedure.


\vspace{-0.2cm}



\subsection{Comparative Analysis of BDCA and Chance-Constrained Solvers}


In light of the superior performance of \gls{BDCA} compared to \gls{DCA}, this section focuses exclusively on its comparison with solvers originating from chance-constrained programming. 
The literature review in Subsection~\ref{sec:chance_constraint_programming} emphasized that a wide range of approaches exists in this area.
In this study, we limit ourselves to algorithms that fulfill the following criteria: (i) compatibility with the sample-based framework, (ii) rigorous convergence theory, (iii) availability of public code implementations, and (iv) use of sequential approximations of the chance constraint rather than conservative convex alternatives. 
Specifically, we selected the \gls{pDCA} recently proposed by \citet{Wang2025}, and the \gls{SCA} derived by \citet{Hong2011}.
Both approaches guarantee convergence to a \gls{KKT} point, and MATLAB code for both algorithms is available via the public GitHub repository of \citet{Wang2025}. 
We translated their implementation into Python to ensure consistency across our experiments.
Note that while \gls{BDCA} and \gls{pDCA} approximate the quantile-based formulation of the chance constraint, \gls{SCA} targets the probability-based formulation using the indicator function. 
Similar to the previous subsection, 250 random starting points are generated for each initialization strategy, and the experiments focus on stricter constraint settings with $r_{\text{min}} \in \{0.96, 0.966, 0.968, 0.97, 0.972 \}$. 
Moreover, runs are terminated after 240 seconds if convergence is not achieved to limit computation costs. 

\begin{figure}[H]%
    \centering
    \includegraphics[width=1\textwidth]{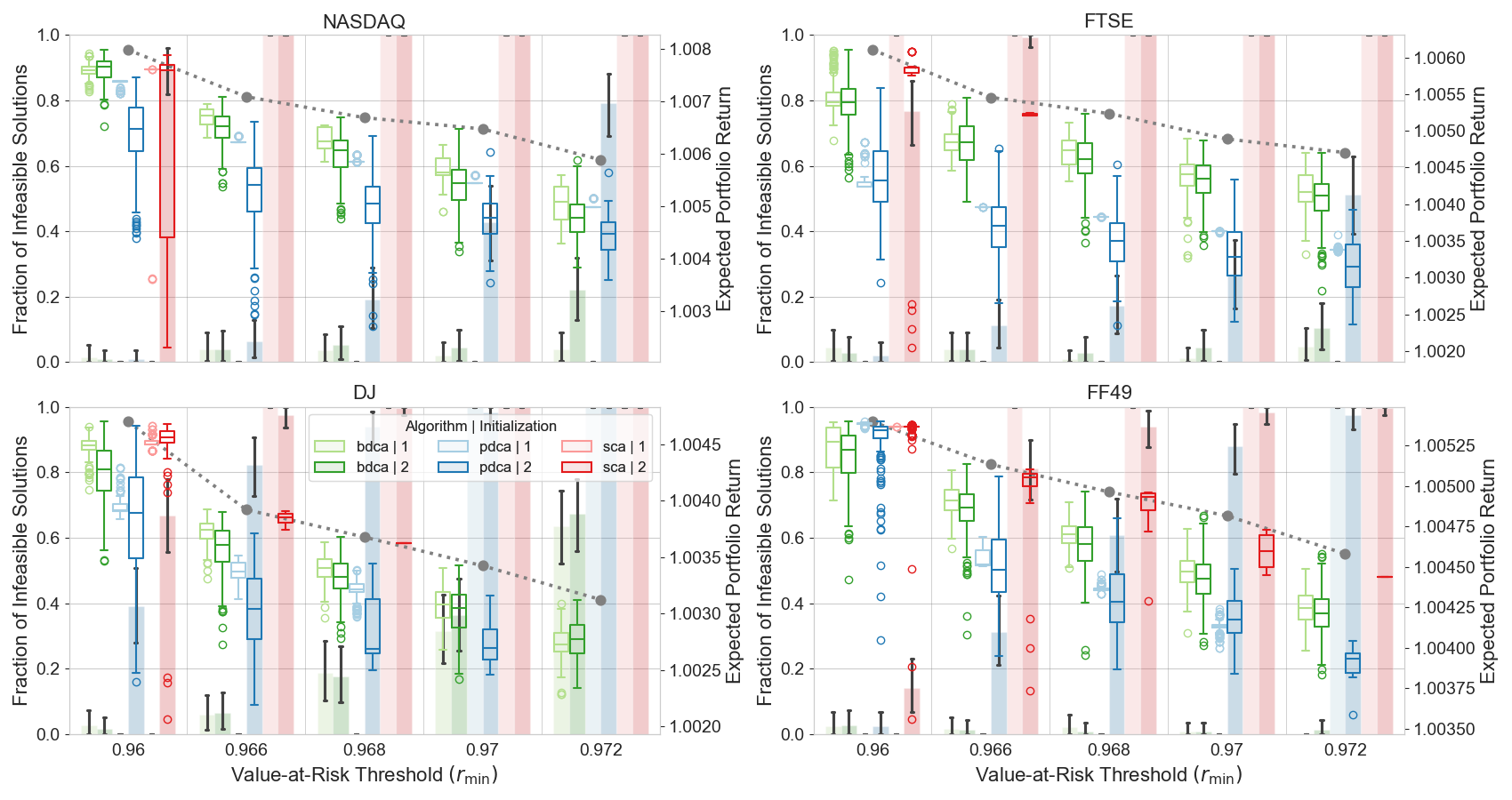}
    \caption{Performance comparison between \gls{BDCA} and chance-constrained solvers using financial datasets and focusing on algorithmic reliability and output quality. 
    Each subplot displays the fraction of infeasible solutions (left y-axis, bar plots) and the expected portfolio returns of all feasible solutions (right y-axis, box plots) across varying $\text{\gls{VaR}}_{X, 0.05}$ levels on the x-axis. 
    Results are grouped by solver type (\gls{BDCA}, \gls{pDCA}, \gls{SCA}) and initialization strategy (1 or 2). 
    The gray horizontal line marks the highest portfolio profit achieved among all feasible solutions, independent of the solver.}
    \label{fig:comp_chance_const_1}%
\end{figure}

Figure~\ref{fig:comp_chance_const_1} compares the fraction of infeasible solutions (left y-axis, bar plots) and the expected portfolio returns of all feasible outcomes (right y-axis, box plots) across the selected settings. 
Among the tested solvers, \gls{SCA} exhibits the strongest dependence on the starting point, as it frequently fails to find feasible solutions for nearly all $r_{\text{min}}$ values and both initialization schemes. 
However, most of its few feasible outcomes tend to be very close to the best-known expected portfolio return for the respective constraint. 
On the other hand, the \gls{pDCA} shows a notable sensitivity dependent on the initialization scheme.
Particularly under the second strategy, it more often leads to infeasible outcomes and the variability of its returns is amplified. 
Independent of the scheme, the median returns of \gls{pDCA} consistently exhibit the largest distance to the best-known feasible solution. 
In contrast, \gls{BDCA} demonstrates the most robust performance across the initialization strategies.
The algorithm consistently identifies feasible solutions even under the most challenging $\text{\gls{VaR}}_{X, 0.05}$ conditions.
While the expected returns of \gls{BDCA} solutions vary to some extent, they are generally more stable and higher than those of \gls{pDCA}.


\vspace{-0.2cm}

\begin{figure}[H]%
    \centering
    \includegraphics[width=1\textwidth]{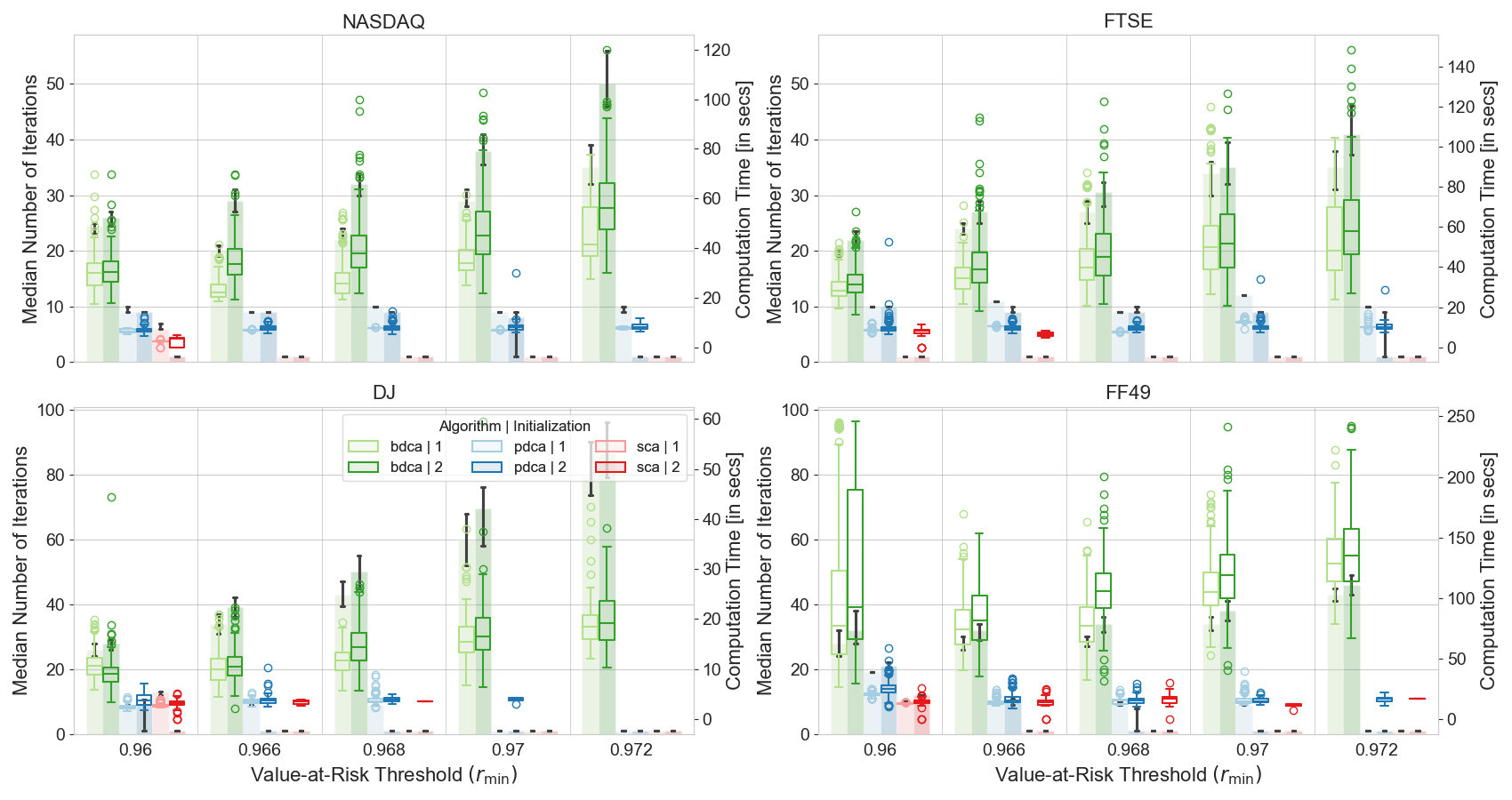}
    \caption{Performance comparison between \gls{BDCA} and chance-constrained solvers using financial data sets and focusing on algorithmic efficiency. 
    Each subplot displays the median number of iterations (left y-axis, bar plots) and the computation time of all feasible solutions (right y-axis, box plots) across varying $\text{\gls{VaR}}_{X, 0.05}$ levels on the x-axis. 
    Results are grouped by solver type (\gls{BDCA}, \gls{pDCA}, \gls{SCA}) and initialization strategy (1 or 2).}%
    \label{fig:comp_chance_const_2}%
\end{figure}

Figure~\ref{fig:comp_chance_const_2} compares the computational effort, measured by the median iteration count (left y-axis, bar plots) and the computation time of all feasible solutions (right y-axis, box plots), across the selected settings. 
Overall, both benchmark solvers require fewer iterations than the \gls{BDCA}.
This then also systematically results in lower computation times for all their feasible outcomes. 
However, especially for the \gls{SCA}, this advantage is offset by its frequent failure to fulfill the different $\text{\gls{VaR}}_{X, 0.05}$ constraints, as discussed previously. 
This imbalance also partially holds for the \gls{pDCA}, as feasibility is only consistently ensured under the first initialization scheme and for the easiest $r_{\text{min}}$ values. 
For more demanding  $\text{\gls{VaR}}_{X, 0.05}$ levels, the \gls{pDCA} provides a better algorithmic efficiency than the \gls{BDCA}, but predominantly offers worse results regarding the output quality and reliability. 
For the \gls{BDCA}, the increased number of iterations and longer computation times are a consequence of the sequential penalty approach. 
Since the initial value for the penalty parameter is deliberately chosen such that the  $\text{\gls{VaR}}_{X, 0.05}$ is not enforced at the beginning, the algorithm naturally requires more iterations until a sufficiently large value is reached. 
Moreover, the subproblem becomes more computationally expensive due to the penalized constraint in the objective function, which introduces a higher number of additional variables.
If it is known in advance that less challenging constraint settings are of interest, the algorithm could alternatively be implemented with the \gls{DC} constraint in non-penalized form to improve the efficiency trade-off. 
Overall, the proposed version of \gls{BDCA} shows a disadvantage in terms of algorithmic efficiency compared to the other approaches, but across performance criteria, it still emerges as the most reliable and robust choice in problems with demanding $\text{\gls{VaR}}_{X, 0.05}$ settings.


\vspace{-0.2cm}


\section{Conclusions}


For over 70 years the \gls{MV} portfolio selection problem has been a topic of great interest in the financial academic literature. 
Although many extensions have been proposed to model the specific requirements of financial markets more accurate, the accompanying increase in complexity often leads to non-convex optimization problems being much harder to solve. 


In this paper, we proposed the \gls{BDCA} for approximately solving a scenario-based portfolio selection problem with a \gls{VaR} constraint. 
For this purpose, we first modified the \gls{DC} representation of the discrete \gls{VaR} presented in \citet{wozabal2010} such that both components are strongly convex.
This modification was then used in a corresponding \gls{DC} program with linear equality and inequality constraints, and where the \gls{VaR} constraint is incorporated via an exact penalty approach.
As in this program both \gls{DC} components are non-smooth, we proposed a novel line search scheme for the \gls{BDCA} to ensure finite termination of the backtracking procedure.  
Afterwards, we proved that the new version of the \gls{BDCA} linearly converges to a \gls{KKT} point for the optimization problem at hand, and outlined how this result can be generalized to polyhedral \gls{DC} problems with linear equality and inequality constraints. 


In the practical part of the paper, we conducted extensive numerical
experiments to assess the performance of the \gls{BDCA}, \gls{DCA} and other solvers from chance-constrained programming for portfolio selection problems with challenging \gls{VaR} constraints. 
By dividing the evaluation into reliability, efficiency, and quality of algorithmic output, we were able to show that the \gls{BDCA} provides the most robust and reliable performance across different initialization strategies and constraint settings.
The algorithm consistently identified the highest number of feasible solutions and had a moderate variability of the expected returns even under the most challenging conditions, while the other approaches failed in these settings.
Due to the open availability of all Python code, we also provided a transparent setup that can be used as starting point for other researchers. 
 

Future work may investigate the incorporation of additional real-world features into the proposed \gls{DC} framework, such as sparsity constraint or transaction costs.
Moreover, although the \gls{BDCA} was more robust than the other approaches, its computation time still scales with the problem dimensions due to the penalty approach. 
Similarly, the median expected returns indicated that the line search procedure increases the solution quality, but also revealed that the algorithm does not necessarily find the global solution. 
Thus, future research could also address the computational efficiency as well as the algorithmic robustness regarding local optima. 
This could further increase the capability of the algorithm to handle challenging data settings and provide the basis for ``.. a way to do it [even] better ... [Thomas A. Edison].



\section{Acknowledgment}

The authors acknowledge the use of the IRIDIS High Performance Computing Facility, and associated support services at the University of Southampton, in the completion of this work.

Further, the authors would like to thank Selin Damla Ahipa{\c{s}}ao{\u{g}}lu, Francisco J. Arag{\'o}n Artacho, Sophie Potts, Tim Marshall-Cox and Samuel Ward for helpful remarks and insightful discussions throughout the development of this work. 

The authors would also like to express their gratitude to three anonymous referees for their constructive comments, which greatly enhanced the quality of this paper.


\addcontentsline{toc}{section}{References}

\setcitestyle{numbers}
\renewcommand{\bibnumfmt}[1]{#1.}

\bibliographystyle{dcu}
\setlength{\bibsep}{4.25pt}
\bibliography{02_Chapters/07_references}


\appendix

\section{Further Numerical Results}\label{app:numerical_results}


\subsection{Sensitivity Analysis of BDCA Input Parameters}


\begin{figure}[H]
    \centering
    \includegraphics[width=1\textwidth]{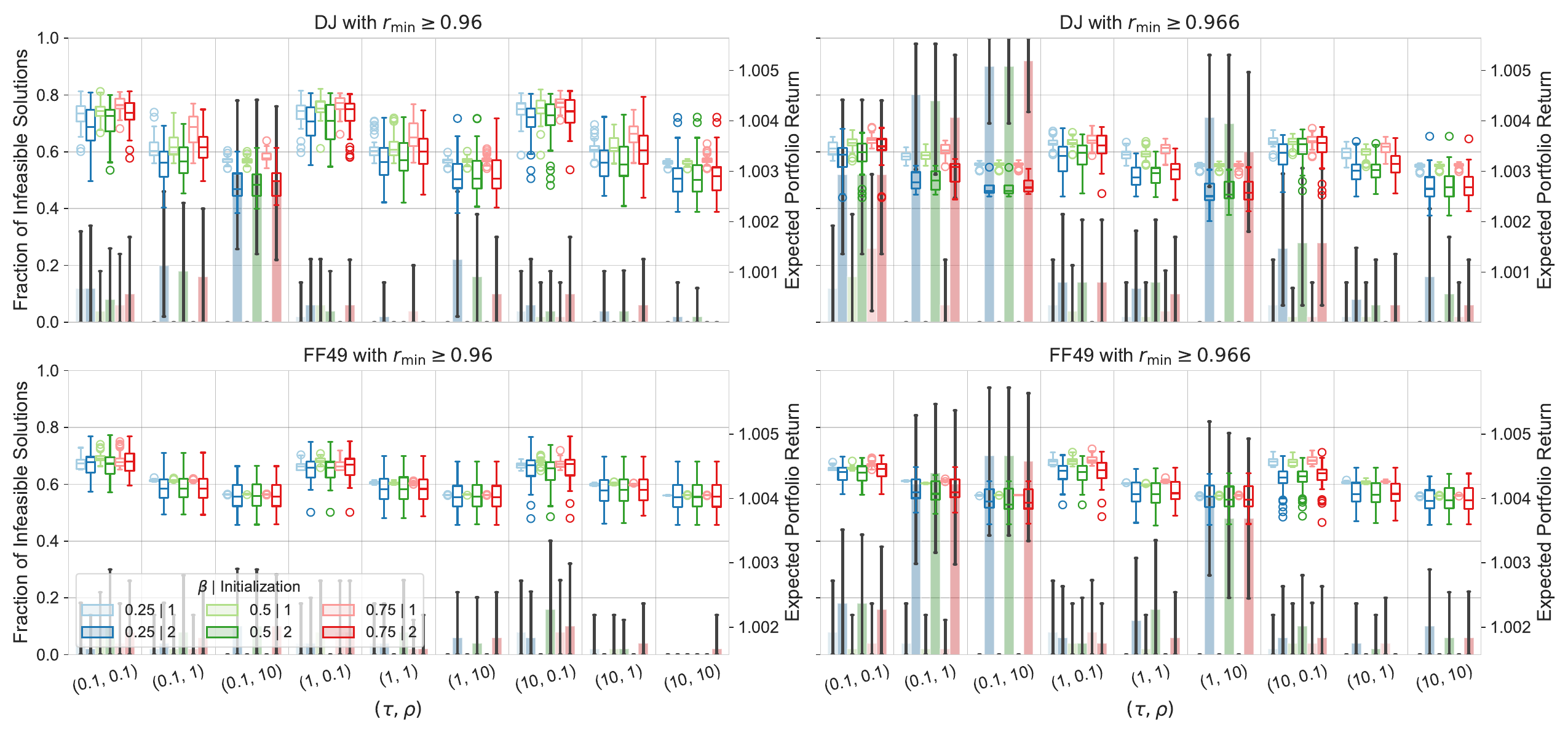}
    \caption{Impact of \gls{BDCA} parameters on feasibility and expected return for the \gls{DJ} and \gls{FF49} data. 
    Each subplot displays the fraction of infeasible solutions (left y-axis, bar plots) and the expected portfolio returns of all feasible outcomes (right y-axis, box plots) across varying combinations of $\tau$ and $\rho$ on the x-axis. 
    Results are grouped by $\beta \in \{0.25, 0.5, 0.75\}$ and initialization strategy (1 or 2). }%
    \label{fig:sensitivity_analysis_bdca_dj_ff49}
\end{figure}


\begin{figure}[H]
    \centering
    \includegraphics[width=1\textwidth]{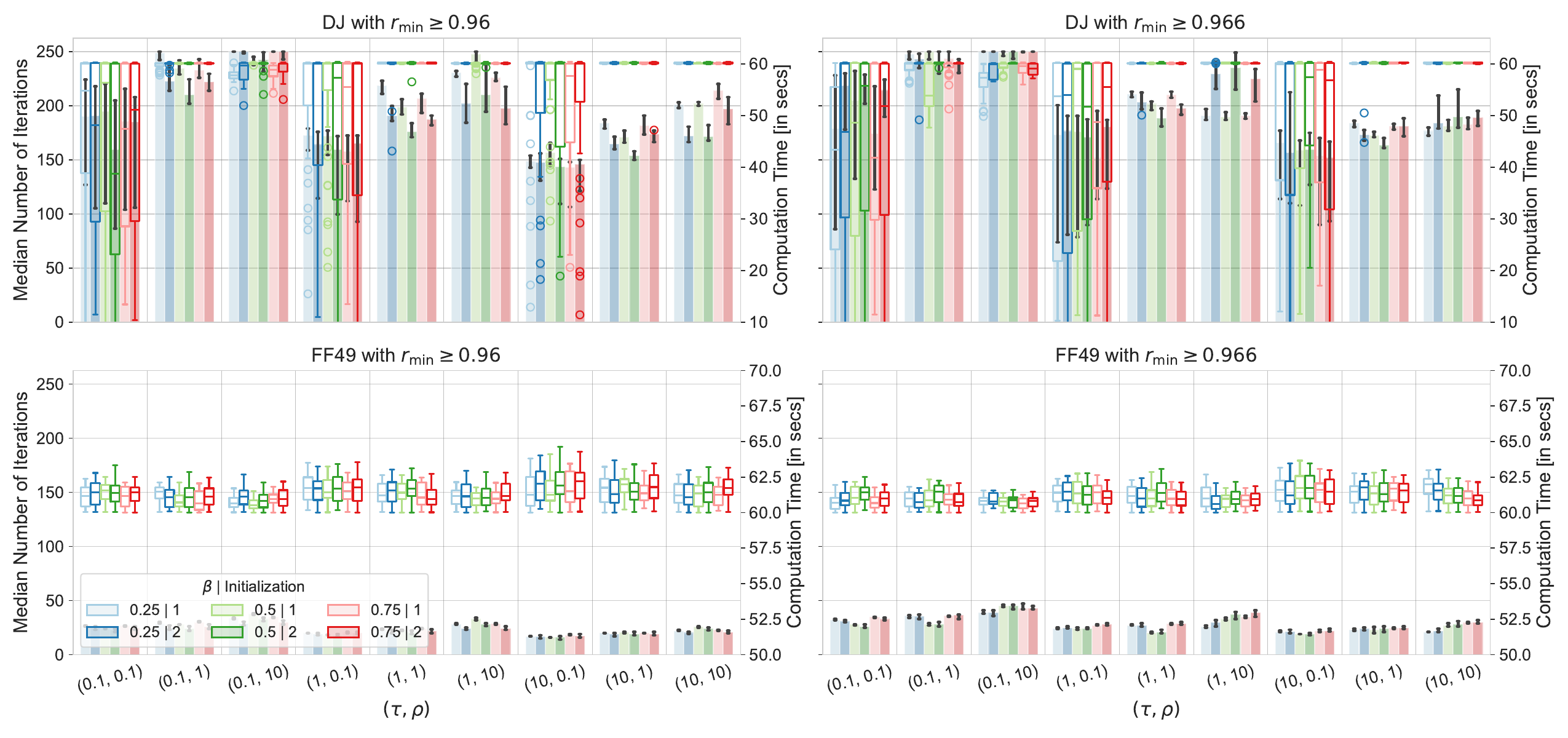}
    \caption{Impact of \gls{BDCA} parameters on the computation time and iteration count for the \gls{DJ} and \gls{FF49} data.
    Each subplot displays the median number of iterations (left y-axis, bar plots) and the computation time of all feasible outcomes (right y-axis, box plots) across varying combinations of $\tau$ and $\rho$ on the x-axis. 
    Results are grouped by $\beta \in \{0.25, 0.5, 0.75\}$ and initialization strategy (1 or 2).}%
    \label{fig:sensitivity_analysis_bdca_2_dj_ff49}
\end{figure}


\subsection{Comparative Analysis of BDCA and DCA}


\begin{figure}[H]%
    \centering
    \includegraphics[scale=0.32]{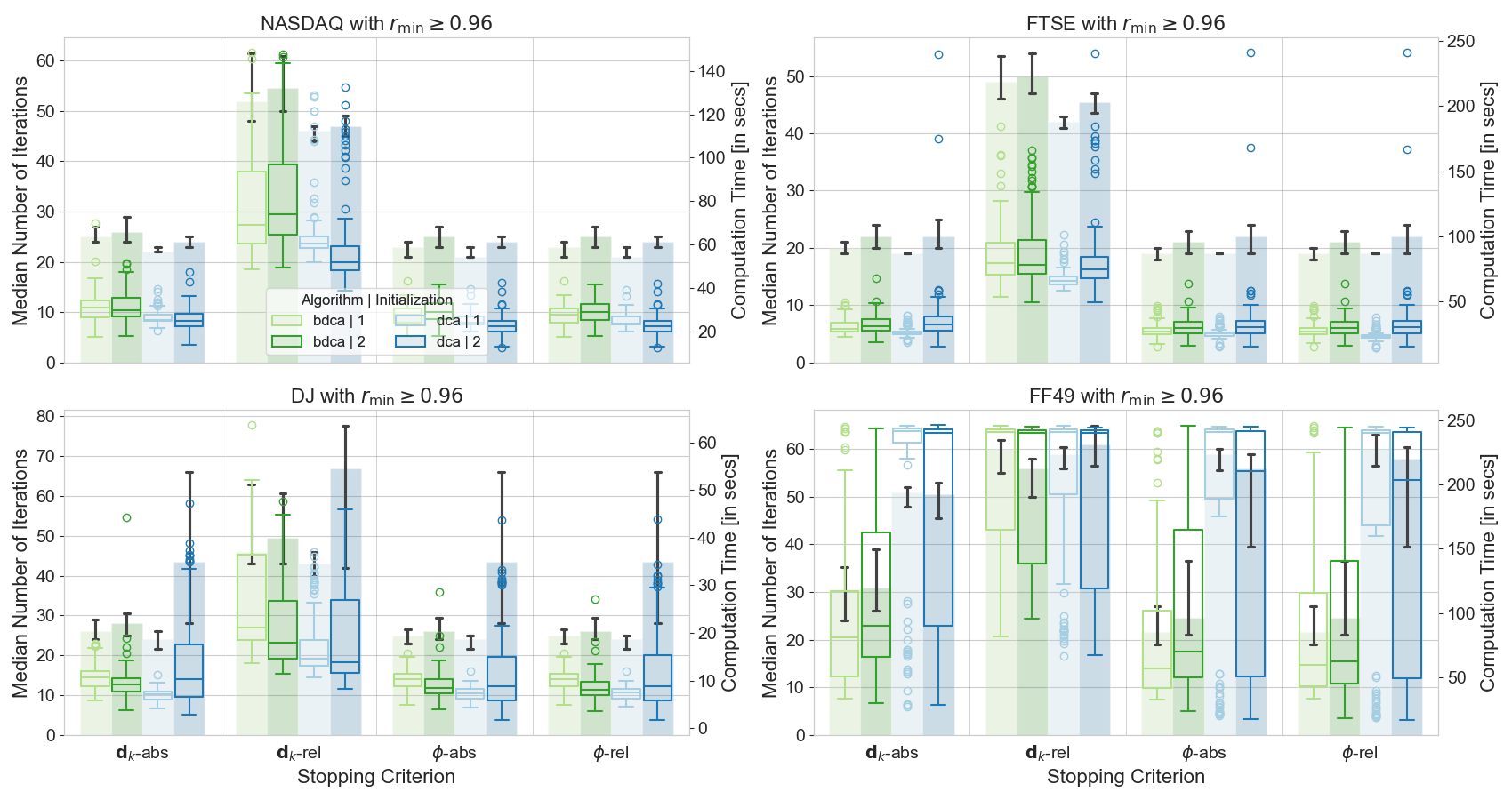}
    \caption{Performance comparison between \gls{BDCA} and \gls{DCA} using financial datasets with $r_{\text{min}} \geq 0.96$ and focusing on algorithmic efficiency depending on different stopping criteria. Each subplot displays the median number of iterations (left y-axis, bar plots) and the computation time of all feasible solutions (right y-axis, box plots) across the different stopping conditions on the x-axis. Results are grouped by solver type (\gls{BDCA}, \gls{DCA}) and initialization strategy (1 or 2).}%
    \label{fig:stop_crit_ct_96}%
\end{figure}


\begin{figure}[H]%
    \centering
    \includegraphics[scale=0.32]{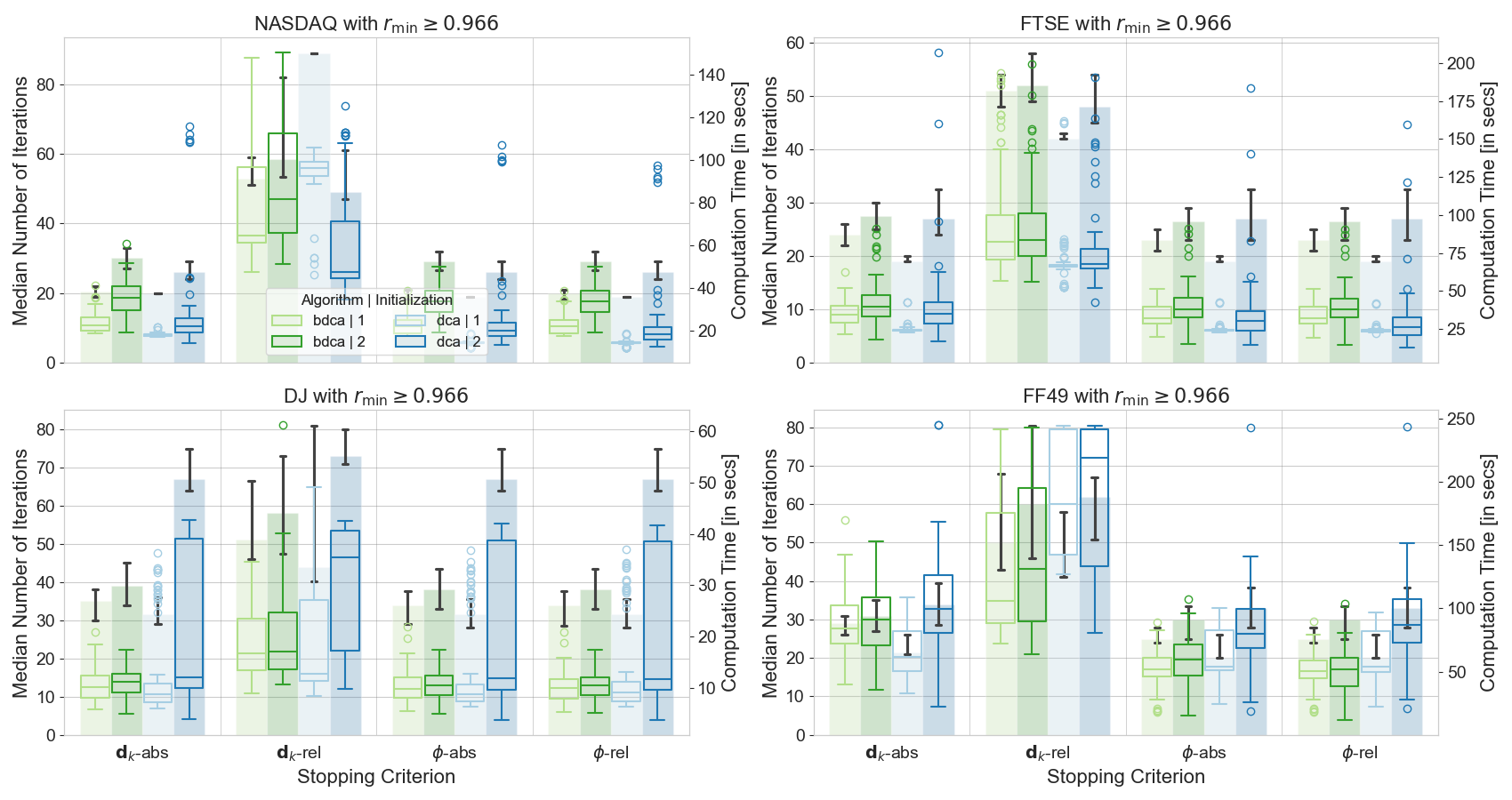}
    \caption{Performance comparison between \gls{BDCA} and \gls{DCA} using financial datasets with $r_{\text{min}} \geq 0.966$ and focusing on algorithmic efficiency depending on different stopping criteria. Each subplot displays the median number of iterations (left y-axis, bar plots) and the computation time of all feasible solutions (right y-axis, box plots) across the different stopping conditions on the x-axis. Results are grouped by solver type (\gls{BDCA}, \gls{DCA}) and initialization strategy (1 or 2).}%
    \label{fig:stop_crit_ct_966}%
\end{figure}


\begin{figure}[H]%
    \centering
    \includegraphics[scale=0.32]{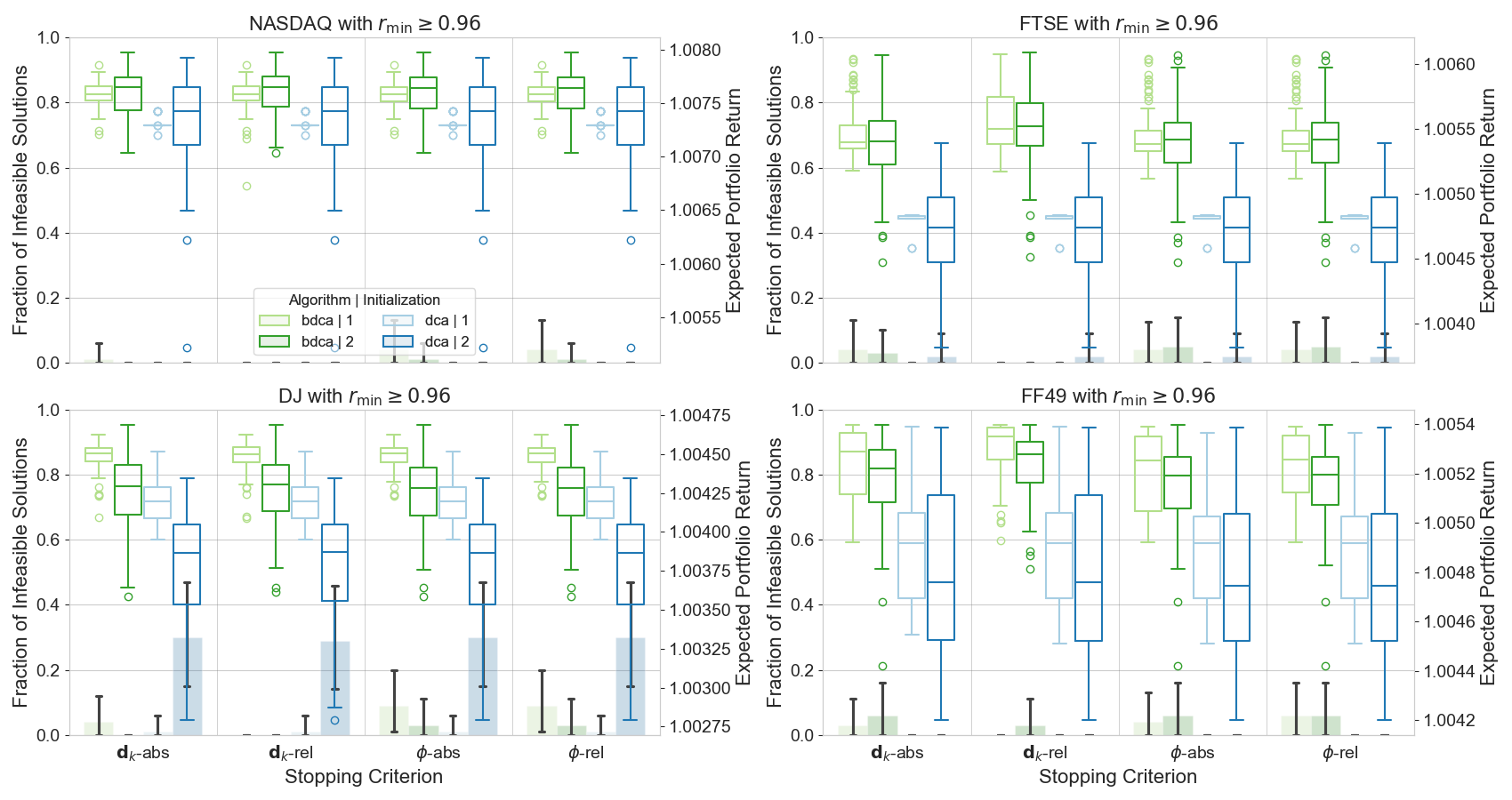}
    \caption{Performance comparison between \gls{BDCA} and \gls{DCA} using financial data sets with $r_{\text{min}} \geq 0.96$ and focusing on algorithmic reliability and output quality depending on different stopping criteria. Each subplot displays the fraction of infeasible solutions (left y-axis, bar plots) and the expected portfolio returns of all feasible outcomes (right y-axis, box plots) across the different stopping conditions on the x-axis. Results are grouped by solver type (\gls{BDCA}, \gls{DCA}) and initialization strategy (1 or 2).}%
    \label{fig:stop_crit_return_96}%
\end{figure}


\begin{figure}[H]%
    \centering
    \includegraphics[scale=0.32]{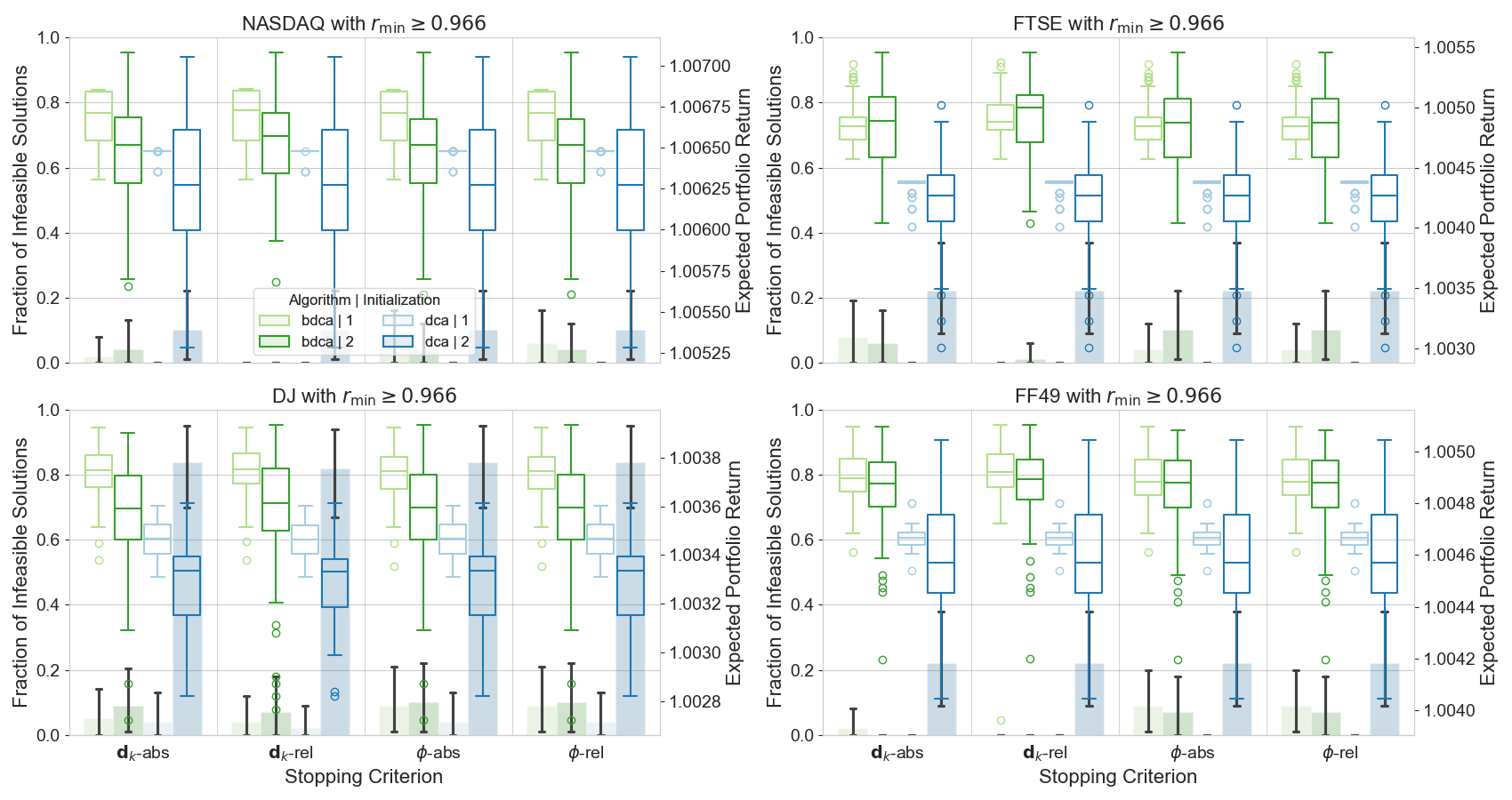}
    \caption{Performance comparison between \gls{BDCA} and \gls{DCA} using financial data sets with $r_{\text{min}} \geq 0.966$ and focusing on algorithmic reliability and output quality depending on different stopping criteria. Each subplot displays the fraction of infeasible solutions (left y-axis, bar plots) and the expected portfolio returns of all feasible outcomes (right y-axis, box plots) across the different stopping conditions on the x-axis. Results are grouped by solver type (\gls{BDCA}, \gls{DCA}) and initialization strategy (1 or 2).}%
    \label{fig:stop_crit_return_966}%
\end{figure}


\newpage
\printglossary[type=\acronymtype, title=List of Abbreviations]
\printglossary

\end{document}